\documentclass[a4paper]{article}
\usepackage{lmodern}
\usepackage[T1]{fontenc}
\usepackage{hyperref}
\usepackage{amsfonts,amsmath,amssymb,amsopn,amsthm}
\usepackage[all]{xy}
\usepackage{vmargin}
\usepackage{makeidx}
\usepackage{graphicx}
\usepackage{tikz}
\usepackage{comment}
\usepackage{titlesec}

\usepackage{color}

\title{Twisted calculus}
\author{Bernard Le Stum \& Adolfo Quir\'os\thanks{Supported by grant MTM2012-35849 from Ministerio de Econom\'{\i}a y Competitividad (Spain).}}
\date{Version of \today}
\newtheorem{thm}{Theorem}[section]
\newtheorem{prop}[thm] {Proposition}
\newtheorem{cor}[thm] {Corollary}
\newtheorem{lem}[thm] {Lemma}
\theoremstyle{definition}%%%%NEEDS \usepackage{amsthm} REMOVE IT IF YOU DO NOT LIKE THIS STYLE
\newtheorem{dfn}[thm] {Definition}

\newenvironment{xmp}[1][Example]{\begin{trivlist} \item[\hskip \labelsep {\bfseries #1}]}{\end{trivlist}}
 \newenvironment{xmps}[1][Examples]{\begin{trivlist} \item[\hskip \labelsep {\bfseries #1}]}{\end{trivlist}}
\newenvironment{pf}[1][Proof]{\begin{trivlist} \item[\hskip \labelsep {\bfseries #1}]}{\end{trivlist}}
\newenvironment{rmk}[1][Remark]{\begin{trivlist} \item[\hskip \labelsep {\bfseries #1}]}{\end{trivlist}}
 \newenvironment{rmks}[1][Remarks]{\begin{trivlist} \item[\hskip \labelsep {\bfseries #1}]}{\end{trivlist}}
\titleformat{\subsection}[runin]{\normalfont\normalsize\bfseries}{\thesubsection}{1em}{}
\numberwithin{equation}{section}
\begin{document}

%\begin{comment}

\maketitle

\bigskip

\begin{center}
\textbf{Abstract}
\end{center}

\medskip
A twisted ring is a ring endowed with a family of endomorphisms satisfying certain relations.
One may then consider the notions of twisted module and twisted differential module.
We study them and show that, under some general hypothesis, the categories of twisted modules and integrable twisted differential modules are equivalent.
As particular cases, one recovers classical results from the theory of finite difference equations or $q$-difference equations. 
\medskip

\tableofcontents

\section*{Introduction}

Finite difference equations (see for example \cite{LevyLessman92}) have been used for a long time to approximate solutions of differential equations.
There  exists also a multiplicative variant called $q$-difference equations (see \cite{DiVizioRamisSauloyZhang03} for example).
Both notions, finite difference equations and $q$-difference equations, can be studied in analogy with the theory of differential equations.
It is possible to give a unified treatment to these different theories through non commutative calculus.
We want to mention in particular Yves Andr\'e's article \cite{Andre01} and the work of Valery Lunts and Alexander Rosenberg (see \cite{LuntsRosenberg97} for example).
Andr\'e's theory introduces the notion of non commutative connections and leads to a beautiful Galois theory.
Lunts and Rosenberg are able to define rings of differential operators in non commutative geometry.
These rings appear naturally in representation theory.

Our approach is based on the quantum philosophy: we want to see the theories of finite difference equations and $q$-difference equations as perturbations of the usual theory of differential equations.
In order to do that, we work over a ring $A$ endowed with a family of endomorphisms satisfying some fixed conditions.
For example, if $A$ is a ring of functions in one variable $x$, and we use one endomorphism $\sigma(x) = x + h$ (and no conditions), we recover the theory of finite difference equations; for $\sigma(x) = qx$, we recover the theory of $q$-difference equations; and, finally, when $\sigma$ is the identity, we recover the theory of usual differential equations.
Even these examples are actually more general than their classical counterpart and encompass the case of positive characteristic as well as the case when $q$ is a root of unity (which is in fact our main concern).

We start with the notion of $E$-twisted ring (resp. module) which is a ring $A$ (resp. an $A$-module $M$) endowed with a family $\underline \sigma := \{\sigma_{i}\}_{i \in E}$ of ring endomorphisms (resp. semi-linear maps) satisfying some conditions such as commutativity, invertibility, existence of roots and so on.
One first shows that such a ring (resp. module) may be seen as a $G$-ring (resp. $G$-$A$-module), which is a ring endowed with an action of a monoid $G$ by ring  endomorphisms (resp. semi-linear maps). In our situation, $G=G(E)$ will be a monoid naturally associated to $E$ and its conditions.

We can do better and introduce the (non commutative) twisted polynomial ring $A[E]_{\underline \sigma}$ associated to $A$ and $E$ (this is the crossed product of $A$ by $G(E)$).
Then an $E$-twisted $A$-module is nothing but an $A[E]_{\underline \sigma}$-module.
For example, if $A = \mathbb C[x]$ and $\sigma(x) = qx$, giving a $\sigma$-twisted module is equivalent to giving a module over the non-commutative polynomial ring $\mathbb C[x, y]$ with the commutation rule $yx = qxy$.
Our aim with the twisted philosophy is to replace non commutative objects with twisted objects that seem easier to handle in practice.

Fix now a commutative base ring $R$. Given a commutative $E$-twisted $R$-algebra $A$, one defines a $\sigma_i$-derivation of $A$ (resp. of an $A$-module $M$) as an $R$-linear map $D_i$ satisfying the twisted Leibnitz rule
\begin{displaymath} 
D_{i}(xs)= D_{i}(x)s + \sigma_{i}(x)D_{i}(s).
\end{displaymath} 
A twisted derivation of $A$ (resp. of an $A$-module $M$) is then a finite sum  of such $D_i$ for some of our endomorphisms $\sigma_i$.

Twisted calculus concerns the study of $A$-modules $M$ endowed with an $R$-linear action of the twisted derivations of $A$ by twisted derivations of $M$.
In order to do that, one can consider the ring of small (or naive) twisted differential operators  $\overline {\mathrm D}_{\underline \sigma}$, which is the smallest ring that contains functions and twisted derivations.
Unfortunately, as  is already the case in the untwisted situation, the category of $\overline {\mathrm D}_{\underline \sigma}$-modules will be too small in general ($p$-curvature phenomenon).

To go further, it is necessary to attach to each endomorphism $\sigma_i$ a specific $\sigma_i$-derivation $D_i$ of $A$: this is what we call a twisted differential algebra.
We may also require some commutation properties that we call the twisted Schwarz conditions.
One can then consider $A$-modules $M$ endowed with an action by $\sigma_i$-derivations of those specific $\sigma_i$-derivations $D_i$.
We also introduce a notion of twisted Weyl algebra as a filtered analog of the twisted polynomial ring $A[E]_{\underline \sigma}$ introduced above, which is a graded ring.

Let us concentrate from now on (with a harmless but convenient change of notation) on the case $E = \underline T := \{T_{1}, \ldots, T_{n}\}$ endowed with the commutation conditions $T_{i}T_{j} = T_{j}T_{i}$.
Then, one can show that there exists a one to one correspondence between twisted differential algebras $(A, \underline \sigma, \underline D)$ that satisfy twisted Schwarz conditions and twisted Weyl algebras $\mathrm D_{\underline \sigma, \underline D}$.
Moreover,  as one of our main results shows, the $\mathrm D_{\underline \sigma, \underline D}$-modules correspond exactly to the $A$-modules with an action of the $D_i$'s by $\sigma_i$-derivations that are integrable (another commutation condition).

We want to apply these results when the $\sigma_i$-derivations $D_i$ span all twisted derivations of $A$.
This is the case when there exists what we call twisted coordinates $x_{1}, \ldots, x_{n}$, for example, when $A$ is a twisted localization of the polynomial ring.
Then, there exists partial $\underline \sigma$-derivations $\partial_{1}, \ldots, \partial_{n}$ and we may consider the twisted Weyl algebra $\mathrm D_{\underline \sigma, \underline \partial}$.
Moreover, one can build a canonical map $A[\underline T]_{\underline \sigma} \to \mathrm D_{\underline \sigma, \underline \partial}$ which is an isomorphism when $\underline \sigma$ is essentially non trivial.
In this case, we obtain an equivalence between $A$-modules endowed with an $R$-linear action of twisted derivations that are integrable, on one hand, and twisted $A$-modules, on the other hand, bringing us back to where we started.

Even if some of our results may be found, usually in a specific context, elsewhere in the mathematical literature, we hope that our abstract and very general description of this twisted calculus will be useful to other mathematicians.
On our side, we plan to apply it to the confluence question: we will show that, under some general hypothesis, all twisted Weyl algebras are topologically isomorphic and, as a consequence, we will prove an equivalence between differential systems and $q$-difference systems.
We will also recast the quantum Simpson correspondence introduced by the first author and M. Gros in \cite{GrosLeStum13} in this new context, by the introduction of twisted divided powers.
Ultimately, one also may expect to understand in some cases the classical Simpson's correspondence as a confluence phenomenon from the quantum case.

Many thanks to Michel Gros for all our fruitful conversations during the preparation of this article.

Throughout the paper, $R$ will denote a commutative ring and our rings will always have a unit.
By an \emph{$R$-algebra}, we mean an $R$-module $A$ endowed with an $R$-bilinear multiplication that turns $A$ into a ring.
This is equivalent to giving a \emph{central} homomorphism $R \to A$.
Note however that homomorphisms of rings that are \emph{not} central also naturally appear in the theory.

%%%%%%%%%%%%%%%%%%%%%%%%%%%%%%%%
%%%%%%%%%%%%%%%%%%%%%%%%%%%%%%%%
\section{Twisted modules}

In this section, we review and extend somehow the formalism of semi-linear algebra.
In order to do that, we need first to introduce or recall some notations and vocabulary.

A \emph{word} of length $n$ on a set $E$, whose elements are usually called \emph{letters}, is a finite sequence $u:=i_{0}i_{1} \cdots i_{n-1}$ of letters of $E$.
More precisely, the set of all words on $E$ is
\begin{displaymath}
W(E) := \coprod_{n \in \mathbb N} \mathrm{Hom}(\{0, \ldots, n-1\}, E).
\end{displaymath}
By composition, any map $f : E \to E'$ will induce a map $W(f) : W(E) \to W(E')$ and this is clearly functorial.
We will denote by $1$ the empty word.
We will also need the notion of \emph{opposite} to a word $u= i_{0}i_{1} \cdots i_{n-1}$: this is the word $u^{op} := i_{n-1}i_{n-2} \cdots i_{0}$.
Alternatively, it is obtained by the symmetry on $W(E)$ which is induced by all the permutations $k \mapsto n-k-1$ of the sets $\{0, \ldots, n-1\}$.

A \emph{(word) condition} on the set $E$ is an ordered pair $(u, v) \in W(E) \times W(E)$ that we will denote by a formal equality ``$u=v$''
(we will not write the quotes when we believe that there is no ambiguity).
A condition of the form ``$ij=ji$'' is called a \emph{commuting condition}, a condition of the form ``$ij=1$'' an \emph{inverting condition} and a condition of the form ``$i^n=j^m$'' is called a \emph{root} condition.
Finally, the \emph{opposite} to the condition ``$u=v$'' is the condition $``u^{\mathrm{op}} = v^{\mathrm{op}}$'' obtained from the corresponding symmetry of $W(E) \times W(E)$.

By definition, a \emph{set of (word) conditions} on $E$ is a subset $\Gamma \subset W(E) \times W(E)$ (the graph of a relation  on $W(E)$).
If we are given two sets $E$ and $E'$, both endowed with a set of conditions $\Gamma$ and $\Gamma'$ respectively, then an \emph{equivariant map} $f : E \to E'$ is a map that preserves conditions, that is,
\begin{displaymath}
(W(f) \times W(f))(\Gamma) \subset \Gamma'.
\end{displaymath}
Pairs $(E, \Gamma)$ made of a set $E$ and a set of word conditions on $E$ make a category with equivariant maps as morphisms.
When the choice of $\Gamma$ is clear from the context, we will simply write $E$ and call $\Gamma$ the \emph{set of $E$-conditions}.
Note that the forgetful functor $(E, \Gamma) \mapsto E$ has an adjoint $E \mapsto E^{\mathrm{triv}} := (E, \emptyset)$ (and also a coadjoint, by the way).
In general, when $\Gamma = \emptyset$, we say that $E$ is \emph{free from conditions}.

We will denote by $\Gamma^{\mathrm{op}}$ the set of opposite conditions to those of $\Gamma$ and $E^{\mathrm{op}} := (E, \Gamma^{\mathrm{op}})$.
We will denote by $\Gamma^{\mathrm{ab}}$ the union of $\Gamma$ and all the commutation conditions and write also $E^{\mathrm{ab}} := (E, \Gamma^{\mathrm{ab}})$.
We let $E^{\pm}$ be the disjoint union $E \coprod E^{\mathrm{op}}$ with the conditions of $\Gamma$, the conditions of $\Gamma^{\mathrm{op}}$ and all the inverting conditions $ii^{\mathrm{op}} = i^{\mathrm{op}}i = 1$ (where we denote by $i^{\mathrm{op}}$ the element $i$ seen as an element of $E^\mathrm{op}$).
Finally, if $S$ is a set of positive integers, let us write $i^{1/n}$ for the element $(i, n)\in E \times S$. Then we let $E^{1/S} := E \times S$ endowed with the root conditions 
\begin{displaymath}
(i^{1/n})^m = (i^{1/n'})^{m'} \quad \mathrm{whenever} \quad \frac mn = \frac {m'}{n'},
\end{displaymath}
plus all the conditions induced by $\Gamma$ on each $E^{1/n} \subset E^{1/S}$ (more precisely, we require $ i_{0}^{1/n} \cdots i_{k-1}^{1/n} =  j_{0}^{1/n} \cdots j_{l-1}^{1/n}$
whenever $i_{0} \cdots i_{k-1} = j_{0} \cdots j_{j-1}$).

When there already exists a multiplication on $E$, which is usually denoted by using juxtaposition, the word notation might be confusing.
This is why we will sometimes apply these considerations to a set of variables $\underline T := \{T_{i}\}_{i \in E}$ rather than to the set $E$ itself.
This is just a matter of notations.

%%%%%%%%%%%%%
\begin{dfn}
Let $E$ be a set endowed with a set $\Gamma$ of (word) conditions.
Then, an \emph{$E$-twisted ring} is a ring $A$ endowed with a family of  ring endomorphisms $\underline \sigma_{A} = \{\sigma_{A,i}\}_{i \in E}$ satisfying
\begin{equation} \label{compos}
\sigma_{A,i_{0}} \circ \cdots \circ \sigma_{A,i_{n-1}} = \sigma_{A,j_{0}} \circ \cdots \circ \sigma_{A,j_{m-1}}
\end{equation}
for all the conditions $i_{0}\cdots i_{n-1} = j_{0} \cdots j_{m-1}$ of $\Gamma$.
It is said to be \emph{inversive} if all $\sigma_{A,i}$'s are bijective.
A \emph{morphism of $E$-twisted rings} $(A, \underline \sigma_{A}) \to (B, \underline \sigma_{B})$ is a ring homomorphism $\varphi : A \to B$ such that
\begin{displaymath}
\forall i \in E, \quad \sigma_{B,i} \circ \varphi = \varphi \circ \sigma_{A,i}.
\end{displaymath}
\end{dfn}
 
 It will be convenient to write $\sigma_{A,u} := \sigma_{A,i_{0}} \circ \cdots \circ \sigma_{A,i_{n-1}}$ when $u = i_{0}\cdots i_{n-1}$ so that the required property becomes $\sigma_{A,u} = \sigma_{A,v}$ for all the conditions $u=v$ of $\Gamma$.
 We also implicitly set $\sigma_{A,1} := \mathrm{Id}_{A}$.
 
Of course, $E$-twisted rings form a category.
It is important to emphasize the fact that this definition depends on the set of conditions $\Gamma$ and we should actually say \emph{$(E, \Gamma)$-twisted ring}.
However, when we don't need to specify $E$, we might simply say \emph{twisted ring}.
It means that $A$ is endowed with a set of endomorphisms satisfying some particular conditions.

Although a twisted ring is a pair $(A, \underline \sigma_{A})$, we will often only use the first letter $A$.
We will also drop the index $A$ and simply write $\sigma_{i}$ for an endomorphism when no confusion can arise.

Before giving some examples of twisted rings, we need to introduce some more vocabulary.

%%%%%%%%%%%%%%%%%%%%%
\begin{dfn}
\begin{enumerate}
\item A \emph{trivially $E$-twisted ring} is an $E$-twisted ring $A$ with $\sigma_{A,i} = \mathrm{Id}_{A}$ for all $i \in E$.
\item An  \emph{$n$-twisted ring} is an $E^{\mathrm{ab}}$-twisted ring where $E := \{1, \ldots, n\}$ and $\Gamma = \emptyset$.
\end{enumerate}
\end{dfn}

%%%%%%%%%%%%%%%%%%%%%
\begin{xmps}
\begin{enumerate}
\item 
Whatever $(E, \Gamma)$ is, the category of trivially $E$-twisted rings is equivalent to the category of usual rings.
This allows us to see the category of rings as a full subcategory of the category of $E$-twisted rings.
\item A $1$-twisted ring is a pair $(A, \sigma_{A})$ where $A$ is a ring and $\sigma_{A}$ is a ring endomorphism of $A$: one can use the one element set $\{T\}$ with no condition at all.
Adding the condition $T^2 = 1$ (resp. $T^2 = T$) for example corresponds to the requirement that $\sigma_{A}$ is a symmetry (resp. a projector). 
\item An $n$-twisted ring is a ring $A$ endowed with $n$ \emph{commuting} ring endomorphisms $\sigma_{A,1}, \ldots, \sigma_{A,n}$: we may use the family $\{T_{1}, \ldots, T_{n}\}$ with all the commutation conditions $T_{i}T_{j}=T_{j}T_{i}$.
More generally, when $E$ is free from conditions, then an $E^{\mathrm{ab}}$-ring is a ring endowed with a family $\{\sigma_{A,i}\}_{i \in E}$ of commuting endomorphisms.
\item An inversive $1$-twisted ring is a pair $(A, \sigma_{A})$ where $A$ is a ring and $\sigma_{A}$ is a ring \emph{automorphism} of $A$.
The full subcategory of inversive $1$-twisted rings is equivalent the category of $T^{\pm}$-twisted ring.
More precisely, $T^{\pm} := \{T, T^{-1}\}$ with the relations $TT^{-1}=T^{-1}T=1$ and we actually  find triples $(A, \sigma_{A}, \sigma_{A}^{-1})$.
More generally, the category of $E^{\pm}$-twisted rings is always equivalent to the category of inversive $E$-twisted rings.
\item If $S$ is a set of positive integers, then $T^{1/S} := \{T^{\frac 1n}\}_{n \in S}$ is endowed with the root conditions
\begin{displaymath}
{(T^{\frac 1n})}^m = {(T^{\frac 1{n'}})}^{m'} \quad \mathrm{whenever} \quad \frac mn = \frac {m'}{n'}.
\end{displaymath}
If $1 \in S$, then a $T^{1/S}$-twisted ring is essentially a $1$-twisted ring endowed with a systems of $n$-th roots $\sigma_{A,n} $ of $\sigma_{A}$ for $n \in S$ (in the sense that $\sigma_{A,n}^m$ only depends on $r := \frac mn \in \mathbb Q$ and $\sigma_{A,n}^n=\sigma_{A}$; see \cite{LeStumQuiros15}).
Again, one can extend this construction to families.
\end{enumerate}
\end{xmps}

Recall that we always denote by $R$ a fixed commutative ring.
Then, we can give some more explicit examples:

%%%%%%%%%%%%%%%%%
\begin{xmps}
\begin{enumerate}
\item If $q_{i}, h_{i} \in R$ for $i = 1, \ldots, n$, then the endomorphisms of the $R$-algebra $R[x_{1}, \ldots x_{n}]$ given by
\begin{displaymath}
\sigma_{i}(x_{j}) = \left\{\begin{array} {l} q_{i}x_{i} + h_{i}\ \mathrm{if}\ j = i \\ x_{j}\ \mathrm{otherwise} \end{array}\right.
\end{displaymath}
turns $R[x_{1}, \ldots x_{n}]$ into an $n$-twisted ring.
\item Let $S$ be a non empty set of positive integers and $\{q_{n}\}_{n \in S}$ a \emph{system of roots} in $R$ (see \cite{LeStumQuiros15} again).
We let $\sigma_{n}$ denote the unique endomorphism of the $R$-algebra $R[x]$ such that $\sigma_{n}(x) = q_{n}x$.
By construction, $R[x]$ endowed with the $\sigma_{n}$'s, is a $T^{1/S}$-twisted ring.
Of course, we can mix this construction with the previous example.
\end{enumerate}
\end{xmps}

%%%%%%%%%%%%%%
\begin{rmks}
\begin{enumerate}
\item One may define the notion of $E$-twisted object in \emph{any} category exactly in the same way.
Cases of particular interest are twisted abelian groups, twisted $A$-modules when $A$ is a fixed ring, twisted $R$-algebras and twisted categories (twisted object of the (quasi-) category of all categories).
We will make this more precise when we meet them.
\item In particular, the notions of trivially-twisted and $n$-twisted object also extend to any category.
\item We can identify the category of $1$-twisted objects of $\mathcal C$ with the kernel of the diagram
\begin{displaymath}
\xymatrix{ \mathrm{Mor}(\mathcal C) \ar@<2pt>[rr]^-{\mathrm{dom}} \ar@<-2pt>[rr]_-{\mathrm{cod}} && \mathcal C}
\end{displaymath}
where $\mathrm{Mor}(\mathcal C)$ denotes the category of all morphisms $u : X \to Y$ and the domain and codomain maps send $u$ respectively to $X$ and $Y$.
Moreover, the canonical common section of $\mathrm{dom}$ and $\mathrm{cod}$ given by $X \mapsto \mathrm{Id}_{X}$ allows the identification of a usual object of $\mathcal C$ with the corresponding trivially $1$-twisted object.
\item 
In his book \cite{Kedlaya10}, Kedlaya calls \emph{difference ring} what we call \emph{$1$-twisted ring} but the terminology \emph{inversive} is borrowed from him.
\end{enumerate}
\end{rmks}

When $\sigma$ is an endomorphism of a ring $A$, we will consider the \emph{pull back functor}
\begin{displaymath}
M \mapsto \sigma^*M := A\ {}_{\sigma}\! \otimes_{A}  M
\end{displaymath}
along $\sigma$ on left $A$-modules.
Unless otherwise specified, we only consider \emph{left} $A$-modules.

%%%%%%%%%%%%%
\begin{dfn} \label{emod}
Let $A$ be an $E$-twisted ring.
Then, an \emph{$E$-twisted $A$-module} is an  $A$-module $M$ endowed with a family of $A$-linear maps $\{\sigma_{M,i} : \sigma_{A,i}^*M \to M\}_{i\in E}$ satisfying
\begin{gather*}
\sigma_{M,i_{0}} \circ \sigma_{A,i_{0}}^*(\sigma_{M,i_{1}}) \circ \cdots \circ \sigma_{A,i_{0}\cdots i_{n-2}}^{*}(\sigma_{M,i_{n-1}}) \circ \sigma_{A,i_{0}\cdots i_{n-1}}^{*}(\sigma_{M,i_{n}})= \\\sigma_{M,j_{0}} \circ \sigma_{A,j_{0}}^*(\sigma_{M,j_{1}}) \circ \cdots \circ \sigma_{A,j_{0}\cdots j_{m-2}}^{*}(\sigma_{M,j_{m-1}})\circ \sigma_{A,j_{0}\cdots j_{m-1}}^{*}(\sigma_{M,j_{m}})
\end{gather*}
for all the $E$-conditions $i_{0}\cdots i_{n} = j_{0}\cdots j_{m}$.
An $E$-twisted $A$-module is said to be \emph{inversive} if all $\sigma_{M,i}$'s are bijective.
A \emph{morphism of $E$-twisted $A$-modules} $(M, \underline \sigma_{M}) \to (N, \underline \sigma_{N})$ is an $A$-linear map $u : M \to N$ such that
\begin{displaymath}
\forall i \in E, \quad \sigma_{N,i} \circ \sigma_{A,i}^*(u) = u \circ \sigma_{M,i}.
\end{displaymath}
\end{dfn}

Of course, $E$-twisted $A$-modules form a category that will be shown later (corollary \ref{abcat}) to be abelian with sufficiently many injective and projective objects.

%%%%%%%%%%%%%%%%%%%%%%%
\begin{xmps}
\begin{enumerate}
\item If $A$ is an $E$-twisted ring, then $A$ itself, together with the identity maps (and \emph{not} with the $\sigma_{i}$'s), is an $E$-twisted $A$-module.
This is even an inversive $E$-twisted $A$-module.
\item In the case $E = \{1, \ldots, n\}$ with only the commutation conditions, we will say \emph{$n$-twisted} again.
Thus, if $A$ is an $n$-twisted ring, an $n$-twisted $A$-module is an $A$-module $M$ endowed with $n$ $A$-linear maps $\sigma_{i} : A\ {}_{\sigma_{A,i}}\! \otimes_{A}  M \to M$ satisfying
\begin{displaymath}
\forall i, j\in\{ 1, \ldots, n\}, \quad \sigma_{i} \circ \sigma_{i}^{*}(\sigma_{j}) = \sigma_{j} \circ \sigma_{j}^{*}(\sigma_{i}).
\end{displaymath}
\item  In particular, if $A$ is a $1$-twisted ring, then a $1$-twisted $A$-module is an $A$-module endowed with an $A$-linear map $\sigma_{M} : A\ {}_{\sigma}\! \otimes_{A}  M \to M$.
It is inversive if and only if $\sigma_{M}$ is bijective.
\end{enumerate}
\end{xmps}

%%%%%%%%%%%%%%%%%%%%%%
 \begin{rmks}
\begin{enumerate}
\item If $A$ is trivially $E$-twisted, then an $E$-twisted $A$-module is the same thing as an $E$-twisted object of the category of $A$-modules, but this is not the case in general due to the pull backs by $\sigma_{A,\underline{i}}$.
\item We mentioned above the notion of $E$-twisted category.
For example, when $A$ is an $E$-twisted ring, we may consider the category of all $A$-modules with the inverse image functors $\sigma_{A,i}^*$.
In general, if $(\mathcal C, \underline \sigma_{\mathcal C}^*)$ is an $E$-twisted category, one can mimic the definition \ref{emod} and obtain the notion of $E$-twisted object of the $E$-twisted category $(\mathcal C, \underline\sigma_{\mathcal C}^*)$.
For example, we will consider below the notions of $\underline \sigma_{A}$-module or $E$-twisted $A$-algebra which also arise from this kind of construction.
\item We can identify the category of $1$-twisted objects of a $1$-twisted category $(\mathcal C, \sigma_{\mathcal C}^*)$ with the kernel of the diagram
\begin{displaymath}
\xymatrix{ \mathrm{Mor}(\mathcal C) \ar@<2pt>[rr]^-{\mathrm{dom}} \ar@<-2pt>[rr]_-{\sigma_{\mathcal C}^* \circ \mathrm{cod}} && \mathcal C.}
\end{displaymath}
\end{enumerate}
\end{rmks}

Recall that if $\sigma$ is an endomorphism of a ring $A$, then the above extension functor $M \mapsto \sigma^*M := A\ {}_{\sigma}\! \otimes_{A}  M$ has a right restriction adjoint $M \mapsto {}\sigma_{*}M := {^{\sigma}M}$.
As an abelian group, ${}^\sigma M$ is identical to $M$ but the action of $A$ is obtained by composition with $\sigma$.
An $A$-linear map $M \to {}^{\sigma} M$ is a \emph{$\sigma$-linear map}: an additive map $\sigma_{M} : M \to M$ such that
\begin{displaymath}
\forall x \in A, \forall s \in M, \quad \sigma_{M}(xs) = \sigma(x)\sigma_{M}(s).
\end{displaymath}

%%%%%%%%%%%%%
\begin{dfn}
If $A$ is an $E$-twisted ring, then a \emph{$\underline \sigma_{A}$-module} is a pair $(M, \underline \sigma_{M})$ where $M$ is an $A$-module and $\underline \sigma_{M} = \{\sigma_{M,i}\}_{i \in E}$ is a family of $\sigma_{A,i}$-linear endomorphisms of $M$ satisfying
\begin{displaymath}
\sigma_{M,i_{0}} \circ \cdots \circ \sigma_{M,i_{n-1}} = \sigma_{M,j_{0}} \circ \cdots \circ \sigma_{M,j_{m-1}}
\end{displaymath}
for all $E$-conditions.
A \emph{morphism of $\underline \sigma_{A}$-modules} $(M, \underline \sigma_{M}) \to (N, \underline \sigma_{N})$ is an $A$-linear map $u : M \to N$ such that
\begin{displaymath}
\forall i \in E, \quad \sigma_{N,i} \circ u = u \circ \sigma_{M,i}.
\end{displaymath}
\end{dfn}

Again, $\underline \sigma_{A}$-modules form a category that we will denote by $\underline \sigma_{A}\mathrm{-Mod}$.

%%%%%%%%%%%%%%%%%%%%%%%%%%%%
\begin{rmk}
If $A$ is an $E$-twisted ring, then we may endow the category dual to that of $A$-modules with the $E$-twisted structure given by the functors $\sigma_{A,i*}$.
We may then consider the category of $E$-twisted objects of this new $E$-category and then take the dual category.
We obtain exactly the category of $\sigma_{A}$-modules.
\end{rmk}

%%%%%%%%%%%%%%%%%%%%%
\begin{prop} \label{nonsens}
If $A$ is an $E$-twisted ring, then the category of $E$-twisted $A$-modules is equivalent (and actually isomorphic) to the category of $\underline \sigma_{A}$-modules.
\end{prop}

\begin{pf}
This is an immediate consequence of adjointness, but let us recall how it works.

If $M$ is an $E$-twisted $A$-module, then the obvious composite map
\begin{displaymath}
M \to A\ {}_{\sigma_{i}}\! \otimes_{A}  M \to M
\end{displaymath}
is $\sigma_{A,i}$-linear;
and conversely, any $\sigma_{A,i}$-linear endomorphism of $M$ extends uniquely to a linear map $A\ {}_{\sigma_{i}}\! \otimes_{A}  M \to M$.
In other words, we have
$$
\sigma_{M,i}(s) = \sigma_{M,i}(1\ {}_{\sigma_{i}}\! \otimes s) \quad \mathrm{and} \quad \sigma_{M,i}(x\ {}_{\sigma_{i}}\! \otimes s) = x\sigma_{M,i}(s).
\quad \Box
$$
\end{pf}

In the future, we might as well call a $\underline \sigma_{A}$-module $M$ an $E$-twisted $A$-module and conversely.
In particular, we will call a $\sigma_{A}$-module \emph{inversive} when the corresponding $E$-twisted $A$-module is.
Be careful that when $M$ is inversive, the semi-linear maps $\sigma_{M,i} : M \to M$ need not be bijective unless $A$ itself is inversive.

%%%%%%%%%%%%%%
\begin{rmks}
\begin{enumerate}
\item If $(\mathcal C, \underline \sigma_{\mathcal C}^*)$ is an $E$-twisted category and we are given a compatible system of adjoint $\sigma_{\mathcal C, i*}$ to the $\sigma_{\mathcal C,i}^*$, then the dual to the category of $E$-twisted objects of $(\mathcal C^{\mathrm{op}}, \underline \sigma_{\mathcal C*})$ is equivalent to the category of $E$-twisted objects of $(\mathcal C, \underline \sigma_{\mathcal C}^*)$.
The proposition is then a particular instance of this equivalence.
\item In his book \cite{Kedlaya10}, Kedlaya calls \emph{(dualizable) difference module} what we call \emph{(inversive) $\sigma_{A}$-module} (or equivalently \emph{(inversive) $1$-twisted module}).
\end{enumerate}
\end{rmks}

If $A \to B$ is a morphism of $E$-twisted rings, then both extension and restriction of scalars preserve $E$-twisted modules.
For extension, one can use the isomorphisms
\begin{displaymath}
B\ {}_{\sigma_{B,i}}\! \otimes_{B}  (B \otimes_{A} M) \simeq B \otimes_{A} (A\ {}_{\sigma_{A,i}}\! \otimes_{A} M),
\end{displaymath}
and for restriction, one may use the canonical morphism $A\ {}_{\sigma_{A,i}}\! \otimes_{A} M \to B\ {}_{\sigma_{B,i}}\! \otimes_{B} M$.
Alternatively, this can be done in the language of $\underline \sigma$-modules.

Also, if $f : E \to E'$ is an equivariant map, then there exists an obvious pull-back functor from $E'$-twisted rings (resp. modules) to $E$-twisted rings (resp. modules).

%%%%%%%%%%%%%%%
\begin{xmps}
\begin{enumerate}
\item The map $\emptyset \to E$ induces a functor that forgets the $E$-structure.
\item If $i \in E$, then the inclusion map $\{i\} \hookrightarrow E$ induces a functor that sends the $E$-twisted ring $(A, \underline \sigma_{A})$ (resp. module $(M, \underline \sigma_{M})$) to a $1$-twisted ring $(A, \sigma_{A,i})$ (resp. module $(M, \sigma_{M,i})$).
\item The map $E \to E^{\pm}$ induces an equivalence between inversive $E$-twisted rings (or modules) and $E^{\pm}$-twisted rings (or modules).
\end{enumerate}
\end{xmps}

As we have already mentioned (and will prove in \ref{abcat}), if $A$ is an $E$-twisted ring, then the category of $\underline \sigma_{A}$-modules is abelian with sufficiently many projective and injective objects and we may therefore define for a $\underline \sigma_{A}$-module $M$,
\begin{displaymath}
\mathrm R\Gamma_{\underline \sigma}(M) := \mathrm{RHom}_{\underline \sigma_{A}\mathrm{-Mod}}(A, M) \quad \mathrm{and} \quad \mathrm H^i_{\underline \sigma}(M) := \mathrm{Ext}^i_{\underline \sigma_{A}\mathrm{-Mod}}(A, M).
\end{displaymath}

%%%%%%%%%%%%%%%%%%%%
\begin{xmps}
\begin{enumerate}
\item We have
\begin{displaymath}
\mathrm H_{\underline \sigma}^0(M) \simeq \cap_{i \in E} \mathrm H_{\sigma_{i}}^0(M) \simeq \{s \in M, \forall i \in E, \sigma_{i}(s) = s\}.
\end{displaymath}
\item If $A$ is a $1$-twisted ring and $M$ is a $1$-twisted $A$-module, one can show that
\begin{displaymath}
\mathrm R\Gamma_{\sigma}(M) \simeq \left[
\xymatrix@R0cm{ M \ar[r]^{1 -\sigma} & M}\right].
\end{displaymath}
This is a complex concentrated in degree $0$ and $1$ whose cohomology is given by
\begin{displaymath}
\mathrm H_{\sigma}^0(M) \simeq \{s \in M, \sigma(s) = s\} \quad \mathrm{and} \quad \mathrm H_{\sigma}^1(M) \simeq M/\{s - \sigma(s), s \in M\}.
\end{displaymath}
\end{enumerate}
\end{xmps}

Note that, if $A$ is an $E$-twisted ring, then $
A^{\underline \sigma = \underline 1} := \mathrm H_{\underline \sigma}^{0}(A)$
is a subring of $A$ called the ring of \emph{$\underline \sigma$-invariants} of $A$.

%%%%%%%%%%%%%%%%
\begin{xmp}
Assume $A = R[x]$ is the polynomial ring in one variable and $\sigma$ is an $R$-endomorphism of $A$.
\begin{enumerate}
\item
If $\sigma(x) = qx$ with $q$ a non trivial primitive $p$-th root of unity in $R$, then $A^{\sigma=1} = R[x^p]$.
\item
If $\mathrm{Char}(R) = p$ and $\sigma(x) = x + 1$, then $A^{\sigma=1} = R[x^p-x]$.
\end{enumerate}
\end{xmp}

Until the end of this section, we will assume (for simplicity) that all the rings are \emph{commutative}.

If $A$ is an $E$-twisted commutative ring, then there always exists an internal tensor product in $\underline \sigma_{A}\mathrm{-Mod}$: if $M$, $N$ are two $\underline \sigma_{A}$-modules, one can endow $M \otimes_{A} N$ with the $\sigma_{i}$-linear endomorphism defined by
\begin{displaymath}
 \sigma_{i}(s \otimes t) = \sigma_{i}(s) \otimes \sigma_{i}(t).
\end{displaymath}
Alternatively, one can use the canonical isomorphism
\begin{equation}
A\ {}_{\sigma_{i}}\! \otimes_{A}  (M \otimes_{A} N) \simeq (A\ {}_{\sigma_{i}}\! \otimes_{A}  M) \otimes_{A} (A\ {}_{\sigma_{i}}\! \otimes_{A} N). \label{tenseq}
\end{equation}
However, in order to turn $\mathrm{Hom}_{A}(M, N)$ into a $\underline \sigma_{A}$-module, we need $M$ to be \emph{inversive} and we may then send $u$ to the composite
\begin{displaymath}
\xymatrix{
\sigma_{i}(u) : M & A\ {}_{\sigma_{i}}\! \otimes_{A}  M \ar[l]_-{\simeq} \ar[r]^{\mathrm{Id} \otimes u} &  A\ {}_{\sigma_{i}}\! \otimes_{A}  N \ar[r] & N
.}
\end{displaymath}
In other words, $\sigma_{i}(u)$ will be characterized by the identity $\sigma_{i}(u) \circ \sigma_{M,i} = \sigma_{N,i} \circ u$.

%%%%%%%%%
\begin{dfn} \label{locadef}
Let $A$ be an $E$-twisted commutative ring.
Then an {$E$-twisted $A$-algebra} (resp. a \emph{$\underline \sigma_{A}$-algebra}) is an $A$-algebra $B$ endowed with a family of $A$-algebra homomorphism $\sigma_{B,i} : A\ {}_{\sigma_{i}}\! \otimes_{A}  B \to B$ (resp. $\sigma_{A,i}$-linear ring endomorphisms $\sigma_{B,i}$) satisfying the $E$-conditions.
\end{dfn}

Again, we obtain two \emph{equivalent} categories (with $A$-linear ring endomorphisms that are compatible with the data).
Actually, giving such a structure is also equivalent to giving a central morphism of $E$-twisted rings $A \to B$.

As a particular case, we will consider the category of $E$-twisted $R$-algebras where $R$ is our trivially twisted commutative base ring.
In other words, an \emph{$E$-twisted $R$-algebra} is an $R$-algebra $A$ endowed with a family of $R$-linear ring endomorphisms $\sigma_{A,i}$ satisfying the $E$-conditions.

It should also be noticed that if $A \to A'$ is any morphism of $E$-twisted commutative rings, then any $E$-twisted $A'$-algebra has a natural structure of $E$-twisted $A$-algebra.
Conversely, if $B$ is an $E$-twisted $A$-algebra, then $\sigma_{A'} \otimes \sigma_{B}$ turns $A' \otimes_{A} B$ into an $E$-twisted $A'$-algebra.

Once more, if $B$ and $C$ are two $E$-twisted $A$-algebras, then $B \otimes_{A} C$ is naturally an $E$-twisted $A$-algebra.

%%%%%%%%%
\begin{dfn} \label{locadef}
Let $A$ be an $E$-twisted commutative ring.
\begin{enumerate}
\item
An $A$-algebra $B$ is an \emph{$E$-twisted quotient} of $A$ if there exists an isomorphism of $A$-algebras $B \simeq A/\mathfrak a$, where $\mathfrak a$ is an ideal in $A$, such that for all $i \in E$, $\sigma_{i}(\mathfrak a) \subset \mathfrak a$.
\item
An $A$-algebra $B$ is an \emph{$E$-twisted localization} of $A$ if there exists an isomorphism of $A$-algebras $B \simeq S^{-1}A$ where $S$ is a multiplicative submonoid of $A$ satisfying for all $i \in E$, $\sigma_{i}(S) \subset S$.
\end{enumerate}
\end{dfn}

%%%%%%%%%%%%%
\begin{prop}
Let $A$ be an $E$-twisted commutative ring and $B$ an $E$-twisted quotient (resp. localization) of $A$.
Then there exists a unique structure of (commutative) $E$-twisted $A$-algebra on $B$.
\end{prop}

\begin{pf}
This follows immediately from the universal property of quotients (resp. localization).
$\quad \Box$
\end{pf}

Conversely, it is clear that if $A \to B$ is a morphism of $E$-twisted commutative rings which is also a quotient (resp. localization) morphism, then $B$ is a $E$-twisted quotient (resp. localization) of $A$.

%%%%%%%%%%%%%%%%%%%%%%%%
\begin{xmps}
\begin{enumerate}
\item If $R$ is an integral domain with fraction field $K$, then any $E$-twisted structure on an $R$-algebra $A$ extends uniquely to $K \otimes_{R} A$.
\item If $R$ is an integral domain with fraction field $K$, and $\sigma$ is a non constant $R$-algebra endomorphism of $R[x]$ (which means that $\sigma(x) \not \in R$), then it extends uniquely to $K(x)$.
\item The $n$-twisted structure given by $\sigma_{i}(x) = q_{i}x$ for some $q_{i} \in R^\times$ on $R[x_{1}, \ldots, x_{n}]$, extends uniquely to $R[x_{1}^{\pm}, \ldots, x_{n}^{\pm}]$.
\item The endomorphism given by $\sigma(x) = x + 1$ of $R[x]$ extends uniquely to $R[x,\frac 1x, \frac 1{x+1}, \ldots,\frac 1 {x+p-1} ]$ when $\mathrm{Char}(R) = p >0$.
\end{enumerate}
\end{xmps}

%%%%%%%%%%%%%%%%%%%%%%%%%%%%%%%%
%%%%%%%%%%%%%%%%%%%%%%%%%%%%%%%%
\section{Twisted polynomial rings}

We will show in this section that twisted modules may always be seen as usual modules over a suitable ring.

Recall that a monoid $G$ is a set endowed with an associative multiplication and a two-sided unit and
that morphisms of monoids preserve all finite products (in particular the unit, which is the empty product).

When $G$ is a monoid, we may endow the set $E_{G} := \{T_{g}\}_{g \in G}$ with the set $\Gamma_{\mathrm{tot}}$ of all conditions $T_{g_{0}}\cdots T_{g_{n-1}} = T_{h_{0}} \cdots T_{h_{m-1}}$ whenever we have an equality $g_{0} \cdots g_{n-1} = h_{0} \cdots h_{m-1}$ in $G$.
One may also use the set $\Gamma_{\mathrm{std}} \subset \Gamma_{\mathrm{tot}}$ of standard conditions
\begin{equation} \label{stdcnd}
T_{g}T_{h}=T_{gh} \quad  \mathrm{for} \quad g, h \in G \quad \mathrm{and} \quad T_{1} = 1.
\end{equation}
Clearly, both maps $G \mapsto (E_{G}, \Gamma_{\mathrm{tot}})$ and $G \mapsto (E_{G}, \Gamma_{\mathrm{std}})$ are functorial in $G$.

Conversely, recall that if $E$ is any set, then $W(E)$ is a monoid for \emph{concatenation}.

The functor $E \mapsto W(E)$ is adjoint to the forgetful functor from monoids to sets.

Now, assume that $E$ is endowed with a set $\Gamma$ of (word) conditions.
This defines a relation $\mathcal R$ on $E$ and we may consider the monoidal equivalence relation  $\bar {\mathcal R}$ generated by $\mathcal R$.
We set
\begin{displaymath}
G(E,\Gamma) := W(E)/(\Gamma) := W(E)/\bar {\mathcal R},
\end{displaymath}
which is a monoid for the quotient structure, and
we obtain a functor $E \mapsto G(E,\Gamma)$ which is adjoint to the above functor $G \mapsto(E_{G}, \Gamma_{\mathrm{tot}})$.

When $\Gamma$ is understood from the context, we will simply write $G(E)$ and denote by $g_{i}$ the class of $i \in E$ in $G(E)$.
However, when we work with a family $\underline T := \{T_{i}\}_{i \in E}$ of variables, we'd rather denote the class of $T_{i}$ by $i$.
The notation we use should in each case be clear from the context.

%%%%%%%%%%%%%%%%
\begin{xmps}
\begin{enumerate}
\item When $E$ is free from conditions, we have $G(E) = W(E)$.
\item As a particular case, we find $G(T) = W(T) \simeq \mathbb N$.
With the condition $T^2=1$, we get $\mathbb Z/2$. And with the condition $T^2=T$, we find $(\{0, 1\}, \times)$.
\item If $T^{\pm}$ denotes the set $\{T, T^{-1}\}$ with the inverting conditions as usual, we see that $G(T^{\pm}) \simeq \mathbb Z$.
More generally, $G(E^{\pm})$ is always the fraction group of the monoid $G(E)$.
\item If we endow $\{T_{1}, \ldots, T_{n}\}$ with the commutation conditions, we find
\begin{displaymath}
G(T_{1}, \ldots, T_{n}) \simeq \mathbb N^n.
\end{displaymath}
More generally, $G(E^{\mathrm{ab}})$ is always the abelian quotient of $G(E)$.
\item Let $S$ be a set of positive integers, then $G(T^{1/S})$ is isomorphic to the submonoid $N := \mathbb N\frac 1 S$ of $\mathbb Q_{\geq 0}$ generated by the set $\{\frac 1n\}_{n \in S}$.
\end{enumerate}
\end{xmps}

%%%%%%%%%%%%%%%%%%%%%
\begin{dfn}
Let $G$ be a monoid, then a \emph{$G$-ring} is a ring $A$ endowed with an action of $G$ by ring endomorphisms.
A \emph{morphism of $G$-rings} is a ring homomorphism $\varphi : A \to B$ which is compatible with the actions.
\end{dfn}

To be more precise, recall that an action of $G$ on $A$ must always satisfy
\begin{displaymath}
\forall g, h \in G, x \in A, \quad (gh).x=g.(h.x) \quad \mathrm{and}\quad  \forall x \in A, \quad 1.x = x,
\end{displaymath}
and we require the extra properties
\begin{displaymath}
\forall g \in G, x, y \in A, \quad g.(x+y)=g.x + g.y \quad \mathrm{and}\quad  g.(xy) = (g.x)(g.y),  \quad \mathrm{and}\quad \forall g  \in G \quad g.1 = 1.
\end{displaymath}
The compatibility condition for homomorphisms means that 
\begin{displaymath}
\forall g \in G, x \in A, \quad \varphi(g.x) = g.\varphi(x).
\end{displaymath}

%%%%%%%%%%%%%%%%%%%%%
\begin{rmks}
\begin{enumerate}
\item Giving an action by ring endomorphisms on the ring $A$ is equivalent to giving a morphism of monoids
\begin{displaymath}
\sigma_{A} : G \to \mathrm{End}_{\mathrm{Rng}}(A).
\end{displaymath}
And the compatibility condition reads
\begin{displaymath}
\forall g \in G, \quad \varphi \circ \sigma_{A}(g) = \sigma_{B}(g) \circ \varphi.
\end{displaymath}
\item One may more generally define a \emph{$G$-object} of a category $\mathcal C$ as an object $X$ endowed with a morphism of monoids $G \to \mathrm{End}(X)$.
Actually, one may see the monoid $G$ itself as a category having exactly one object (with elements of $G$ as endomorphisms), and giving a $G$-object is equivalent to giving a functor $G \to \mathcal C$.
For example, we may consider $G$-sets, $G$-monoids, $G$-modules, $G$-$A$-modules (for the trivial action on $A$) or $G$-categories.
\item The category of $G$-sets is a topos and giving a $G$-ring is equivalent to giving a ring in this topos.
Similarly, giving a $G$-module for example, is equivalent to giving an abelian group in the topos of $G$-sets.
Anyway, as a consequence, we see that all limits (resp. colimits) exist in the category of $G$-rings and, actually, the underlying ring is the limit (resp. colimit) ring.
\item As we noticed before, one may define a $G$-category as a $G$-object of the category of all categories.
For example, the category of $A$-modules has a structure of $G$-category given by pull back when $A$ is a $G$-ring.
And the dual category has such a structure for push-out (right $G$-action).
Another example of $G$-category is given by the slice category $\underline G$ of $G$ (over the only object of $G$ as a category), whose objects are the elements of $G$ and morphisms are all the $h: kh \to k$ (with obvious composition).
One makes $G$ act naturally on the objects ($g.h = gh$) and trivially on the maps ($g.h = h$).
When $G$ is right-cancellative (for example $G = \mathbb N^n$) and endowed with its natural order, one can identify $\underline G$ with the category associated to the ordered set $G$ and the action of $G$ on $\underline G$ is the natural one.
\end{enumerate}
\end{rmks}

%%%%%%%%%%%%%%%%%%%
\begin{prop}
If $E$ is a set with conditions and $G := G(E)$, then the category of $E$-twisted rings is equivalent (and even isomorphic) to the category of $G$-rings.
\end{prop}

\begin{pf}
This is an immediate consequence of the universal property of $G(E)$: any map $E \to \mathrm{End}(A)$ which is compatible with the conditions on $E$ and the total set of conditions on the monoid $\mathrm{End}(A)$ extends uniquely to a morphism of monoids $G(E) \to \mathrm{End}(A)$.
$\quad \Box$
\end{pf}

It follows that the category of $E$-twisted rings only depends on $G(E)$.

%%%%%%%%%%%%%%%%
\begin{xmps}
\begin{enumerate}
\item If $E$ is any set, giving a $W(E)$-ring is equivalent to giving a ring $A$ together with a family $\{\sigma_{A,i}\}_{i \in E}$ of endomorphisms of $A$.
\item Giving an $\mathbb N$-ring (resp. a $\mathbb Z$-ring, a $\mathbb Z/2$-ring, a $\{0, 1\}$-ring) is equivalent to giving a ring $A$ together with a ring endomorphism  (resp. a ring automorphism, a ring symmetry, a ring projector) $\sigma_{A}$.
\item Giving an $\mathbb N^n$-ring (resp. a $\mathbb Z^n$-ring) an is equivalent to giving a ring $A$ together with $n$ commuting endomorphisms (resp. automorphisms) $\sigma_{A,i}$.
\item If $N$ is a submonoid of $\mathbb Q_{>0}$ containing $1$, giving an $N$-ring is equivalent to giving a ring $A$ together with a ring endomorphism $\sigma_{A}$, and a compatible family of $n$-roots $\sigma_{A,n}$ of $\sigma_{A}$ whenever $\frac 1n \in N$.
\end{enumerate}
\end{xmps}

In order to build the twisted polynomial rings, we need to recall a few basic results and set up some notations.
If $A$ is ring and $E$ is any set, we denote by $AE := \oplus_{i\in E}Ai$ the free $A$-module on the basis $E$.
The functor $E \mapsto AE$ is adjoint to the forgetful functor from $A$-modules to sets.
Moreover, there exists an obvious canonical map
\begin{equation} \label{naiv}
\xymatrix@R0cm{
A \times E \ar[r] & AE\\ (x, i) \ar[r] & xi.
}
\end{equation}
Also, when $G$ is a monoid, there exists a natural ring structure on $\mathbb ZG$ and the functor $G \mapsto \mathbb ZG$ is adjoint to the forgetful functor from rings to monoids.

Assume now that $A$ is a $G$-ring for some monoid $G$.
Then, in particular, we can see $A$ as a $G$-monoid (for multiplication) and  consider the semi-direct product $H := A \rtimes_{\sigma} G$ which is a monoid with $A \times G$ as underlying set.
The canonical map $H \to AG$ extends uniquely to a surjective additive homomorphism $\pi : \mathbb ZH \to AG$.

%%%%%%%%%%%%%%%
\begin{lem} \label{tows}
If $A$ is a $G$-ring and $H := A \rtimes_{\sigma} G$, then the kernel of the canonical map $\pi : \mathbb ZH \to AG$ is a two-sided ideal.
\end{lem}

\begin{pf}
By definition, an element 
\begin{displaymath}
\sum_{(x,g)} n_{x,g} (x, g) \in \mathbb ZH
\end{displaymath}
is in the kernel  of $\pi$ if and only if
\begin{displaymath}
\forall g \in G, \quad \sum_{x} n_{x,g}x = 0
\end{displaymath}
in $A$.
Moreover, we have for any $(y, h)  \in A \rtimes_{\sigma} G$,
\begin{displaymath}
\sum n_{x,g} (x, g) (y,h) = \sum n_{x,g} (x\sigma_{g}(y), gh).
\end{displaymath}
In order to show that the kernel is a right ideal, we need to check that
\begin{displaymath}
\forall k \in G, \quad \sum_{\substack{(x,g)\\ gh=k }} n_{x,g} x\sigma_{g}(y) = 0.
\end{displaymath}
But for each $g \in G$ such that $gh=k$, we have $ \sum_{x} n_{x,g}x = 0$ and we are done.
We follow the same process in order to show that the kernel is a left ideal.
If $(y, h)  \in A \rtimes_{\sigma} G$, we have
\begin{displaymath}
(y,h) \sum_{(x,g)} n_{x,g} (x, g)  = \sum_{(x,g)} n_{x,g} (y\sigma_{h}(x), hg).
\end{displaymath}
and we need to show that we always have
\begin{displaymath}
\sum_{\substack{(x,g)\\ gh=k }} n_{x,g} y\sigma_{h}(x) = 0.
\end{displaymath}
But for each $g \in G$ such that $gh=k$, we have
\begin{displaymath}
\sum_{x} n_{x,g}y\sigma_{h}(x) = y\sigma_{h}(\sum_{x} n_{x,g}x) = 0. \quad \Box
\end{displaymath}
\end{pf}

%%%%%%%%%%%%%%%%%%%
\begin{prop} 
If $A$ is a $G$-ring, then there exists a unique bilinear multiplication on the underlying additive abelian group of $AG$ that extends the multiplication of $A$ and $G$ and such that
\begin{displaymath}
\forall a \in A, g \in G, \quad ga = \sigma_{g}(a)g.
\end{displaymath}
This multiplication turns $AG$ into a ring.
\end{prop}

\begin{pf}
This follows from lemma \ref{tows} since this multiplication is induced by the quotient map $\pi : \mathbb ZH \to AG$.
$\quad \Box$
\end{pf}

%%%%%%%%%%%%%%%%%%%
\begin{dfn}
\begin{enumerate}
\item If $A$ is a $G$-ring, then the \emph{crossed product} of $A$ by $G$ is the free $A$-module $AG$ endowed with the ring structure of the proposition.
\item If $A$ is an $E$-twisted ring, then the \emph{$E$-twisted polynomial ring} $A[E]_{\underline \sigma}$ is the crossed product of $A$ by $G(E)$.
\end{enumerate}
\end{dfn}

Actually, we should (and might) write $A[E, \Gamma]_{\underline \sigma}$ and not merely $A[E]_{\underline \sigma}$ because this ring strongly depends on the word conditions.
When $A$ is trivially twisted, we will simply write $A[E]$ or $A[E, \Gamma]$.

%%%%%%%%%%%%%%%%
\begin{rmks}
\begin{enumerate}
\item When the action of $G$ on $A$ is trivial, then the crossed product of $A$ by $G$ is the usual algebra of the monoid $G$ which is usually denoted by $A[G]$.
In order to avoid confusion, we will not use any specific notation to denote the crossed product of $A$ by $G$ in general. 
\item When $\Gamma = \emptyset$,  our conventions say that $A[E]$ denotes the ring of \emph{non-commutative} polynomials on $E$.
On the other hand, the usual polynomial ring, which is generally denoted by $A[E]$, is actually for us $A[E^{\mathrm{ab}}]$.
We might nevertheless use the same notation for both notions when we believe that there is no risk of confusion.
\item In \cite{Sabbah93}, Claude Sabbah uses some kind of twisted polynomial rings to in order to study $q$-difference equations.
More precisely, he uses the ring $A := K[x_{1}, \ldots, x_{n}]$ with $\sigma_{i}(x_{j}) = q^{\partial_{i,j}}x_{j}$ where $K$ a field of characteristic zero and $q$ a non zero element of $K$.
Note that he also considers the case $A = K[x_{1}^{\pm}, \ldots, x_{n}^\pm]$ and $\sigma^\pm_{i}(x_{j}) = q^{\pm\partial_{i,j}}x_{j}$.
\end{enumerate}
\end{rmks}

%%%%%%%%%%%%%%%%
\begin{xmps}
\begin{enumerate}
\item If $A$ is a $1$-twisted ring, then $A[T]_{\sigma}$ is the non-commutative polynomial ring in one variable over $A$ with the commutation rule $Tx=\sigma(x)T$ for $x \in A$ (the Ore extension of $A$ by $\sigma$ and $0$).
\item As a particular case, if $A = R[S]$ is endowed with $\sigma(S) = qS$, then $A[T]_{\sigma}$ is the standard \emph{non-commutative polynomial ring} $R[T, S]_{q}$ over $R$ with the commutation rule $TS = qST$ which is used to defined the quantum plane.
\item If $A$ is an $n$-twisted ring, then $A[T_{1}, \cdots, T_{n}]_{\underline \sigma}$ is the free $A$ module with basis all $T_{1}^{k_{1}} \cdots T_{n}^{k_{n}}$ for $k_{1}, \ldots, k_{n} \in \mathbb N$ endowed with the usual multiplication rules for the variables and $T_{i}x=\sigma_{i}(x)T_{i}$ for $x \in A$ and $i = 1, \ldots, n$.
We call it the \emph{twisted $n$-polynomial ring} on $A$ (note that the commutation of variables is built-in).
\item If $A$ is an inversive $1$-twisted ring, then $A[T^{\pm}]_{\sigma}$ is the non-commutative ring of Laurent polynomials with the commutation rules $Tx=\sigma(x)T$ and $T^{-1}x =\sigma^{-1}(x)T^{-1}$ for $x \in A$.
This extends to several variables and gives rise to \emph{twisted Laurent polynomials}.
\item  Let $S$ be a set of positive integers.
 Then $A[T^{1/S}]_{\sigma}$ is the non commutative ring of Puiseux polynomials with powers in $\mathbb N\frac 1S$, with the commutation rules $T^rx=\sigma_{n}^m(x)T^r$ for $x \in A$ and $r = \frac mn$.
Note that there exists an isomorphism
\begin{displaymath}
A[T^{1/S}]_{\sigma} \simeq \varinjlim_{\{T \to T^n\}_{n\in S}} A[T]_{\sigma_{n}}
\end{displaymath}
This description extends also to several variables and gives rise to \emph{twisted Puiseux polynomials}.
\end{enumerate}
\end{xmps}

The twisted polynomial rings come with the following universal property:

%%%%%%%%%%%%%%%%%%
\begin{prop} \label{unipo}
Let $A$ be an $E$-twisted ring.
Let $\varphi : A \to B$ be a ring homomorphism and $\{y_{i}\}_{i \in E}$ a family of elements of $B$ such that $y_{i_{0}} \ldots y_{i_{n-1}} = y_{j_{0}} \ldots y_{j_{m-1}}$ for any $E$-condition $i_{0}\cdots i_{n} = j_{0} \cdots j_{n}$.
If we have
\begin{displaymath}
\forall  i \in E, \forall x \in A, \quad y_{i}\varphi(x) = \varphi( \sigma_{i}(x))y_{i},
\end{displaymath}
then there exists a unique ring homomorphism $\Phi : A[E]_{\underline \sigma} \to B$ that extends $\varphi$ and such that for all $i \in E$, we have $\Phi(g_{i}) = y_{i}$.
\end{prop}

\begin{pf}
In order to check this assertion, one can use the various universal properties that occur in the construction of and $A[E]_{\underline \sigma}$.
$\quad \Box$
\end{pf}

%%%%%%%%%
\begin{cor}
If $A$ is an $E$-twisted ring, then there exists a canonical isomorphism
\begin{displaymath}
A[E^{\mathrm{triv}}]_{\underline \sigma}/I \simeq A[E]_{\underline \sigma}
\end{displaymath}
where $I$ is the two-sided ideal generated by the conditions of $E$.
\end{cor}

More precisely, $I$ denotes the ideal generated by all the $i_{0}\ldots i_{n-1} - j_{0}\ldots j_{m-1}$ whenever ``$i_{0}\ldots i_{n-1} = j_{0}\ldots j_{m-1}$'' is an $E$-condition.

\begin{pf}
Both rings satisfy the same universal property.
$\quad \Box$
\end{pf}

Note that when $G(E)$ is an abelian monoid, for example when $E = E^{\mathrm{ab}}$, one can use the twisted \emph{commutative} polynomial ring on the left hand side.

%%%%%%%%%%%%%%%%%
\begin{prop} \label{chang1}
If $A \to B$ is a morphism of $E$-twisted rings, then there exists a canonical isomorphism
\begin{displaymath}
B \otimes_{A} A[E]_{\underline \sigma} \simeq B[E]_{\underline \sigma}.
\end{displaymath}
\end{prop}

\begin{pf}
Both rings share the same universal property.
$\quad \Box$
\end{pf}

As an illustration, we can compute the center of $A[T]_{\sigma}$ when $A$ is an integral domain:

%%%%%%%%%%%%%%%
\begin{prop}
Let $A$ be a $1$-twisted integral domain.
Then,
\begin{enumerate}
\item If there exists $p > 0$ such that $\sigma_{A}^p = \mathrm{Id}_{A}$, and $p$ is the smallest such integer, then the centralizer of $A$ in $A[T]_{\sigma}$ is $A[T^p]$ and the center of $A[T]_{\sigma}$ is equal to $A^{\sigma = 1}[T^p]$.
\item Otherwise, the centralizer of $A$ in $A[T]_{\sigma}$ is equal to $A$ and the center of $A[T]_{\sigma}$ is equal to $A^{\sigma = 1}$.
\end{enumerate}
\end{prop}

Of course, with some extra conventions, we can see the second statement as a particular case of the first one (case $p = 0$).

\begin{pf}
In order to lighten the notations, we simply denote by $\sigma$ the endomorphism of $A$.
If $\sum x_{k}T^k \in A[T]_{\sigma}$ commutes with all $x \in A$, we must have for all $k \in \mathbb N$, $x_{k}\sigma^k(x) = x_{k}x$.
Since we assume that $A$ is an integral domain, this can happen only if $x_{k} = 0$ or $\sigma^k = \mathrm{Id}$.
Thus, we see that the centralizer of $A$ in $A[T]_{\sigma}$ is reduced to $A$ unless there exists $k > 0$ such that $\sigma^k = \mathrm{Id}$.
If $p$ is the smallest such integer, then necessarily, $k \equiv 0 \mod p$ and the assertions on the centralizer of $A$ are proved.
Now, an element $\sum x_{k}T^k \in A[T]_{\sigma}$ is in the center of the ring if and only if it is in the centralizer of $A$ and commutes with $T$.
This last condition means that for all $k \in \mathbb N$, $x_{k}$ commutes with $T$ and this is equivalent to $\sigma(x_{k}) = x_{k}$.
In other words, we get that $x_{k} \in A^{\sigma = 1}$.
$\quad \Box$
\end{pf}

%%%%%%%%%%%%%%%%%%%%%
\begin{dfn}
Let $A$ be a $G$-ring.
Then, a \emph{$G$-$A$-module} is an $A$-module $M$ endowed with an action of $G$ by semi-linear maps.
A \emph{morphism of $G$-$A$-modules} is a homomorphism $u : M \to N$ which is compatible with the actions.
\end{dfn}

Semi-linearity means that,  beside the standard properties
\begin{displaymath}
\forall g, h \in G, s \in M, \quad (gh).s=g.(h.s) \quad \mathrm{and}\quad  \forall s \in M, \quad 1.s = s,
\end{displaymath}
we must also have
\begin{displaymath}
\forall g \in G, s, t \in M, \quad g.(s+t)=g.s + g.t \quad \mathrm{and} \quad \forall g \in G, x \in A, s \in M,  g.(xs) = (g.x)(g.s).
\end{displaymath}
The compatibility condition means that 
\begin{displaymath}
\forall g \in G, x \in A, \quad u(g.x) = g.u(x).
\end{displaymath}

%%%%%%%%%%%%%%%%%%%%%
\begin{rmks}
\begin{enumerate}
\item  One can define in general a \emph{$G$-object} of a $G$-category $\mathcal C$ as a $G$-functor $\underline G \to \mathcal C$ (recall that one interprets $G$ as a category and denote by $\underline G$ its slice category).
When $A$ is a $G$-ring and $\mathcal C$ denotes the $G$-category of $A$-modules with all the pull-back functors, this is equivalent to the notion of $G$-$A$-module.
And if we use the push-out and dualize both at the beginning and at the end, we get also an equivalent category.
\item If we see a $G$-ring $A$ as a ring in the topos of $G$-sets, then a $G$-$A$-module is nothing but a module on the $G$-ring $A$.
And similarly for $G$-$A$-algebras for example.
Note that, as a consequence, the category of $G$-$A$-modules is abelian with enough injectives (we will give another proof of this fact below).
Also there exists an internal tensor product and an internal Hom.
Finally, all limits (resp. colimits) exist in the category of $G$-$A$-modules and the underlying $A$-module is the limit (resp. colimit) $A$-module.
\end{enumerate}
\end{rmks}

%%%%%%%%%%%%%%%%%%%
\begin{prop}
If $E$ is a set with conditions and $G := G(E)$, then the category of $E$-twisted $A$-modules is equivalent (and even isomorphic) to the category of $G$-$A$-modules.
\end{prop}

\begin{pf}
Again, this is an immediate consequence of the universal property of $G(E)$.
$\quad \Box$
\end{pf}

In particular, we see that the category of $E$-twisted $A$-modules only depends on $G(E)$.

%%%%%%%%%%%%%%%%%
\begin{rmk}
Under this equivalence, we may see an $E$-twisted ring $A$ (resp. $A$-module $M$) as a ring (resp. an $A$-module) in the topos of $G$-sets.
One can check that the above internal tensor product $M \otimes_{A} N$ corresponds to the usual topos tensor product.
Note that the topos internal hom (that always exists) does not have in general $\mathrm{Hom}_{A}(M, N)$ as underlying $A$-module unless $M$ is strict.
Also, we may interpret $\underline \sigma$-cohomology as topos (or equivalently sheaf) cohomology.
\end{rmk}

%%%%%%%%%%%%%%%
\begin{prop} \label{eqtwp}
Let $A$ be an $E$-twisted ring.
Then the category of $E$-twisted $A$-modules is equivalent (and even isomorphic) to the category of $A[E]_{\underline \sigma}$-modules.
\end{prop}

\begin{pf}
Again, this is completely formal but we can work out the details.
First of all, since both categories only depend on $G := G(E)$, we may assume that $E = \{T_{g}\}_{g \in G}$ with the standard conditions.
Giving an $A$-module structure on an abelian group $M$ is equivalent to giving a morphism of rings $\varphi : A \to \mathrm{End}_{\mathrm{Ab}}(M)$.
Extending the $A$-module structure to an $A[E]_{\underline \sigma}$-module structure is equivalent to extending $\phi$ to a morphism of rings $\Phi : A[E]_{\underline \sigma} \to \mathrm{End}_{\mathrm{Ab}}(M)$.
Thanks to proposition \ref{unipo}, giving $\Phi$ is equivalent to giving a family of $\sigma_{M,g} \in \mathrm{End}_{\mathrm{Ab}}(M)$ for all $g \in G$ satisfying 
\begin{displaymath}
\forall g, h \in G,\quad \sigma_{M,g} \circ \sigma_{M,h} = \sigma_{M, gh} \quad \mathrm{and} \quad \sigma_{M,1} = \mathrm{Id}_{M},
\end{displaymath}
and
\begin{displaymath}
\forall g \in G, x \in A, s \in M, \quad \sigma_{M,g}(xs) = \sigma_{A,g}(x)\sigma_{M,g}(s).
\end{displaymath}
Clearly, this is equivalent to giving a semi-linear action of $G$ on $A$, which in turn is equivalent to an $E$-twisted structure.
$\quad \Box$
\end{pf}
 
 \begin{rmk}
From section 6.3 of Peter Hendriks' thesis \cite{hendriks96}, one obtains a beautiful classification of $1$-twisted modules on $\mathbb C(z)$ when $\sigma(z) = qz$ with $q$ a root of unity.
 \end{rmk}

%%%%%%%%%%%%%%%%%%%%
\begin{cor} \label{abcat}
If $A$ is an $E$-twisted ring, then the category of $E$-twisted $A$-modules is an abelian category with sufficiently many projective and injective objects.$\quad \Box$
\end{cor}

It follows from the proposition that if $M$ is an $E$-twisted $A$-module, then
\begin{equation}
\mathrm R\Gamma_{\underline \sigma}(M) = \mathrm{RHom}_{A[E]_{\underline \sigma}}(A, M). \label{homis}
\end{equation}

\begin{xmp} One may use the complex
\begin{displaymath}
\left[A[T]_{\sigma} \stackrel {1 - T}\longrightarrow A[T]_{\sigma}\right],
\end{displaymath}
which is a free left resolution of $A$, to describe $H^{i}_{\sigma}(M)$ when $M$ is a $1$-twisted $A$-module.
\end{xmp}

One can define, even when $M$ is \emph{not} inversive, the internal hom of two $E$-twisted $A$-modules $M$ and $N$ as
\begin{displaymath}
\mathrm{Hom}_{\underline \sigma\mathrm{-Mod}}(A[E]_{\underline \sigma} \otimes_{A} M, N),
\end{displaymath}
where $A[E]_{\underline \sigma}$ is endowed with its natural action $\sigma_{i}(\sum a_{g}g) = \sum \sigma_{i}(a_{g})g$ and $A[E]_{\underline \sigma} \otimes_{A} M$ is endowed with the internal tensor product.

Finally, note also that if $A \to B$ is any morphism of rings, then the adjoint functors between the categories of $E$-twisted modules on $A$ and $B$ are obtained by extension and restriction through the canonical map $A[E]_{\underline \sigma} \to B[E]_{\underline  \sigma}$.

%%%%%%%%%%%%%%%%%%%%%%%%%%
%%%%%%%%%%%%%%%%%%%%%%%%%%%
\section{Twisted derivations}\label{twder}

In this section, we generalize to the case of twisted algebras the notions of derivation and small (or naive) differential operator.
From now on, we stick to the \emph{commutative} case.

Let $A$ be a commutative $R$-algebra.
If $M$ and $N$ are two $A$-modules, we will denote by $\mathrm{Hom}_{R}(M, N)$ the $R$-module of all $R$-linear maps from $M$ to $N$.
Recall that $\mathrm{Hom}_{R}(M, N)$ has two $A$-module structures coming from the $A$-modules structures of $M$ and $N$ respectively.
In the particular case $M=N$, we will write $\mathrm{End}_{R}(M)$ which is actually an $R$-algebra for composition.
We will implicitly consider multiplication by $x \in A$ as an endomorphism of $M$ (or $N$).
Then, the $A$-modules structures on $\mathrm{Hom}_{R}(M, N)$ are given by composition on the left or composition on the right with these morphisms.
We will also implicitly consider an element $s \in M$ as an $A$-linear map $A \to M$ (the unique $A$-linear map sending $1$ to $s$).

All this applies in particular to the case $M = N = A$ so that $\mathrm{End}_{R}(A)$ will denote the $R$-algebra of all of \emph{$R$-linear} endomorphisms of $A$ (and \emph{not} only $R$-algebra endomorphisms).
We will identify $A$ with an $R$-subalgebra of $\mathrm{End}_{R}(A)$ (sending $x \in A$ to multiplication by $x$ inside $A$).
We will try to avoid any confusion that might arise from this identification.

%%%%%%%%%%%%%%%%
\begin{dfn} Let $A$ be a commutative $R$-algebra, $\sigma$ a ring endomorphism of $A$ and $M$, $N$ two $A$-modules.
If $\varphi \in \mathrm{Hom}_{R}(M, N)$ and $x \in A$, then the \emph{$\sigma$-twisted bracket} of $\varphi$ and $x$ is
\begin{displaymath}
[\varphi, x]_{\sigma} = \varphi \circ x - \sigma(x) \circ \varphi \in  \mathrm{Hom}_{R}(M,N).
\end{displaymath}
\end{dfn}

We might simply say \emph{$\sigma$-bracket} in the future.
Thus, we set
\begin{displaymath}
\forall s \in M,\quad [\varphi, x]_{\sigma}(s) = \varphi(xs) - \sigma (x) \varphi(s).
\end{displaymath}

In the case $\sigma = \mathrm{Id}_{A}$, we will simply write $[\varphi, x]$ so that
\begin{displaymath}
[\varphi, x] = \varphi \circ x - x \circ \varphi
\end{displaymath}
as usual.

%%%%%%%%%%%%%%%
\begin{rmks}
\begin{enumerate}
\item Formally, it would not be necessary to introduce $\sigma$-brackets because the $\sigma$-bracket $[x, \varphi]_{\sigma}$ on $\mathrm{Hom}_{R}(M,N)$ is identical to the usual bracket $[x, \varphi]$ on $\mathrm{Hom}_{R}(M,{}^\sigma N)$.
\item It is sometime convenient to use the vocabulary of bimodules (see section 1.4.2 of \cite{Andre01} or \cite{Bourbaki70}, Chapter 3, section 10):
if $P$ is an $A$-bimodule, the \emph{bracket} of $s \in P$ and $x \in A$ is defined by
\begin{displaymath}
[s, x] = s \cdot x - x \cdot s \in P.
\end{displaymath}
Thus, we see that our $\sigma$-bracket is a particular case of general bimodule bracket applied to the bimodule $P := \mathrm{Hom}_{R}(M,{}^\sigma N)$.
\item
Note also that, for fixed $\sigma$, any $A$-module $M$ may be seen as an $A$-bimodule with right action given by $s \cdot x = xs$ and left action given by $x \cdot s = \sigma(x)s$ (called a \emph{$\sigma$-sesquimodule} in \cite{Andre01}).
In this case, we see that the bimodule bracket on $M$ is given by $[x, s] = ys$, where $y = x - \sigma(x)$ is an element that will play an increasing role later.
\end{enumerate}
\end{rmks}

%%%%%%%%%%%%%%%%%%%%
\begin{lem}
Let $A$ be a commutative $R$-algebra, $\sigma$ a ring endomorphism of $A$ and $M$, $N$ two $A$-modules.
Then, for fixed $x \in A$ the map $\varphi \mapsto [\varphi, x]_{\sigma}$, from $\mathrm{Hom}_{R}(M, N)$ to itself, is $A$-linear both on the left and on the right.
\end{lem}

\begin{pf}
We have for $x \in A$ and $\varphi, \psi \in \mathrm{Hom}_{R}(M, N)$,
\begin{displaymath}
[\varphi + \psi, x]_{\sigma} = \varphi \circ x + \psi \circ x - \sigma (x)  \circ \varphi - \sigma (x)  \circ \psi
\end{displaymath}
\begin{displaymath}
= \varphi \circ x - \sigma (x)  \circ \varphi + \psi \circ x - \sigma (x)  \circ \psi = [\varphi, x]_{\sigma} + [\psi, x]_{\sigma}.
\end{displaymath}
Also, if $x,y \in A$ and $\varphi \in \mathrm{Hom}_{R}(M, N)$, we see that
\begin{displaymath}
[y \circ \varphi, x]_{\sigma} = y \circ \varphi \circ x - \sigma (x) \circ y \circ \varphi
\end{displaymath}
\begin{displaymath}
= y \circ \varphi \circ x - y \circ \sigma (x) \circ \varphi = y \circ [\varphi, x]_{\sigma}
\end{displaymath}
and
\begin{displaymath}
[\varphi  \circ y, x]_{\sigma} = \varphi \circ y \circ x - \sigma (x) \circ \varphi  \circ y
\end{displaymath}
\begin{displaymath}
= \varphi \circ x \circ y - \sigma (x) \circ \varphi  \circ y = [\varphi, x]_{\sigma} \circ y. \quad \Box
\end{displaymath}
\end{pf}

%%%%%%%%%%
\begin{dfn} \label{sigmadev}
Let $A$ be a commutative $R$-algebra and $M$ an $A$-module.
If $\sigma$ is an endomorphism of $A$, then a \emph{$\sigma$-derivation of $A$ into $M$} is an element $D \in \mathrm{Hom}_{R}(A,M)$ such that
\begin{displaymath}
\forall x \in A, \quad [D,x]_{\sigma}  = D(x).
\end{displaymath}
More generally, if $A$ is an $E$-twisted commutative $R$-algebra, then a \emph{$\underline \sigma$-derivation of $A$ into $M$} is a finite sum $D := \sum_{i\in E} D_{i}$ where $D_{i}$ is a $\sigma_{i}$-derivation of $M$.
In the case $M = A$, we call $D$ a \emph{$\underline \sigma$-derivation of $A$ over $R$}.
\end{dfn}

For later use, note that an $R$-linear map $D$ is a $\sigma$-derivation if and only if it satisfies the \emph{$\sigma$-Leibnitz rule}
\begin{equation} \label{sigleib}
\forall x, y \in A, \quad D(xy) = yD(x) + \sigma(x)D(y).
\end{equation}
Note also that if $D$ is a $\underline \sigma$-derivation, then we have $D(a) = 0$ for all $a \in R$.
This follows from the fact that if $D$ is a $\sigma$-derivation, we will have
\begin{displaymath}
D(a) = [D, a]_{\sigma} = 0
\end{displaymath}
since $D$ is $R$-linear.
Finally, note that the notion of $\underline \sigma$-derivation strongly depends on $E$ and not only on $G(E)$ since, in fact, the word conditions do no play any role in the definition.

We can give explicit formulas for derivating powers:

%%%%%%%%%%
\begin{prop} \label{genD1}
Let $\sigma$ be an endomorphism of a commutative $R$-algebra $A$ and $D : A \to M$ a $\sigma$-derivation.
Then,
\begin{enumerate}
\item We have
\begin{displaymath}
\forall k \in \mathbb N, \quad D (x^{k} )= \sum_{j=0}^{k-1} x^j\sigma(x)^{k-1-j}D(x).
\end{displaymath}
\item If we set $y := x - \sigma(x)$, we have also
\begin{displaymath}
\forall k \in \mathbb N, \quad D (x^{k} ) = \sum_{j=0}^{k-1} {k \choose j} (-y)^{k-1-j}x^jD(x).
\end{displaymath}
\end{enumerate}
\end{prop}

\begin{pf}
The first assertion is easily proved by induction on $k$:
\begin{displaymath}
D (x^{k+1}) = x^k D (x) + \sigma(x) D (x^k) =  x^kD(x) + \sigma(x) \sum_{j=0}^{k-1} x^j\sigma(x)^{k-1-j}D(x)
\end{displaymath}
\begin{displaymath}
=  \left(x^k + \sum_{j=0}^{k-1} x^j\sigma(x)^{k-j}\right)D(x) = \sum_{j=0}^{k} x^j\sigma(x)^{k-j}D(x).
\end{displaymath}

The second assertion is also proved by induction on $k$ using $\sigma(x) = x - y$.
We will have
\begin{displaymath}
D (x^{k+1}) = x^k D (x) + \sigma(x) D (x^k) =  x^kD(x) + (x-y)  \sum_{j=0}^{k-1} {k \choose j} (-y)^{k-1-j}x^jD(x)
\end{displaymath}
and we can compute
\begin{displaymath}
 x^k + (x-y)  \sum_{j=0}^{k-1} {k \choose j} (-y)^{k-1-j}x^j
=  x^k +  \sum_{j=0}^{k-1} {k \choose j} (-y)^{k-1-j}x^{j+1} +  \sum_{j=0}^{k-1} {k \choose j} (-y)^{k-j}x^j
\end{displaymath}
\begin{displaymath}
=  x^k +  \sum_{j=1}^{k} {k \choose j-1} (-y)^{k-j}x^{j} +  \sum_{j=0}^{k-1} {k \choose j} (-y)^{k-j}x^j
\end{displaymath}
\begin{displaymath}
=  x^k + kx^k + (-y)^k + \sum_{j=1}^{k-1} {k+1 \choose j}(-y)^{k-j}x^j = \sum_{j=0}^{k} {k+1 \choose j} (-y)^{k-j}x^j. \quad \Box
\end{displaymath}
\end{pf}

%%%%%%%%%%%%%%%%%%%%%%
 
\begin{xmps}
\begin{enumerate}
\item If $\sigma$ is any endomorphism of the polynomial $R$-algebra $A := R[x]$, then there exists a unique $\sigma$-derivation $\partial_{\sigma}$ of $A/R$ such that $\partial_{\sigma}(x) = 1$ (see proposition \ref{genD1} below).
It is given by
\begin{equation} \label{sigform}
\partial_{\sigma} (x^{n} ) = \sum_{i+j+1=n} x^i\sigma(x)^{j}.
\end{equation}
\item Assume that $A := R[x_{1}, \dots, x_{n}]$ is endowed with $n$ endomorphism $\sigma_{i}$ such that $\sigma_{i}(x_{j}) = x_{j}$ whenever $i \neq j$.
Then there exists for each $i = 1, \ldots, n$, a unique $\sigma_{i}$-derivation $\partial_{i}$ of $A/R$ such that
\begin{equation} \label{sigformn}
\forall j = 1, \ldots, n, \quad \partial_{i}(x_{j}) = \left\{\begin{array} {cl} 1 & \mathrm{if}\ j=i \\ 0 & \mathrm{otherwise} \end{array}\right.
\end{equation}
Actually, this is the unique $\sigma_{i}$-derivation of $A/R_{i}$ satisfying $\partial_{i}(x_{i}) = 1$ if we let $R_{i} := R[x_{1}, \ldots, \widehat {x_{i}}, \ldots, x_{n}]$.
\item If $K$ is a field and $\sigma$ is field endomorphism of $K(x)/K$, then there exists a unique $\sigma$-derivation on $K(x)/K$ with $\partial_{\sigma}(x) = 1$.
It is given by
\begin{equation}\label{sigquot}
\partial_{\sigma}\left(\frac yz\right) = \frac {\partial_{\sigma}(y)z - y\partial_{\sigma}(z)}{z\sigma(z)}
\end{equation}
for $y, z \in K[x]$, $K$-linearity and formula \eqref{sigform} for positive powers of $x$.
This case extends to several variables as well.
\item \label{xp4}
For this example, we recall that we introduced in \cite{LeStumQuiros15}, when $q \in R$ and $m \in \mathbb N$, the \emph{quantum  integer}
\begin{displaymath}
(m)_{q} := 1 + q + \cdots q^{m-1}.
\end{displaymath}
Now, assume that $A := R[x]$ and let $\sigma(x) = qx$, then the unique $\sigma$-derivation on $A$ with $\partial_{\sigma}(x) = 1$ is given by
\begin{equation} \label{formit}
\partial_{\sigma}(x^n) = (n)_{q}x^{n-1} \quad \mathrm{for}\ n \in \mathbb N.
\end{equation}
Again, this extends to more variables.
\item 
When $q \in R^\times$ and $m \in \mathbb N$, we may also consider the \emph{quantum integer}
\begin{displaymath}
(-m)_{q} := -\frac 1q - \frac 1{q^2} - \cdots - \frac 1{q^{m}}.
\end{displaymath}
If A := $R[x,x^{-1}]$ and we let $\sigma(x) = qx$, then there exists again a unique $\sigma$-derivation on $A$ with $\partial_{\sigma}(x) = 1$.
It is given by formula \eqref{formit} with $n \in \mathbb Z$.
And this works as well with more variables.
\item If $A := R[x]$ and we let $\sigma(x) = x + h$ with $h \in R$, then the unique $\sigma$-derivation with $\partial_{\sigma}(x) = 1$ is given by
\begin{displaymath}
\partial_{\sigma} (x^{n} ) = \sum_{i+j+1=n} {n \choose i}h^{j}x^i.
\end{displaymath}
\end{enumerate}
\end{xmps}

The next lemma gives a generic example:

%%%%%%%%%%%%%%%%%%
\begin{lem} \label{easy}
If $\sigma$ is an endomorphism of a commutative $R$-algebra $A$, then the operator $1 - \sigma$ is a $\sigma$-derivation of $A/R$.
 \end{lem}

\begin{pf}
We have for all $x, y \in A$,
\begin{displaymath}
(1 -\sigma)(x)y + \sigma(x)(1-\sigma)(y) = (x -\sigma(x))y + \sigma(x)(y-\sigma(y))
\end{displaymath}
\begin{displaymath}
= xy - \sigma(x)y + \sigma(x)y - \sigma(xy) = xy -\sigma(xy) =  (1 - \sigma)(xy). \quad \Box
\end{displaymath}
\end{pf}

%%%%%%%%%%%%%%%%%%
\begin{rmks}
\begin{enumerate}
\item One may introduce more generally (as explained in section 1.4.2 of \cite{Andre01}) the notion of derivation $D : A \to M$ of a ring into a \emph{bimodule} by requiring that $[D, x] = D(x)$ (we use here the bimodule bracket).
This is of course compatible with our definition in the case of a sesquimodule.
Note that giving a derivation $D : A \to M$ into a bimodule is equivalent to giving a section $x \mapsto (x, D(x))$ of the augmentation map $U^1(M) = A \oplus M \to A$ where $U^1(M)$ denotes the (non commutative) tensor algebra of $M$ modulo $M^{\otimes 2}$.
\item When $A$ is inversive, one might also meet \emph{symmetric $\sigma$-derivations} defined by the condition
\begin{equation} \label{sigleib}
\forall x, y \in A, \quad D(xy) = \sigma(y)D(x) + \sigma^{-1}(x)D(y)
\end{equation}
as in \cite{Hodges90}, p. 160.
For example, if $A = R[x]$ and $\sigma(x) = v^2x$ with $v \in R^\times$, then there will exist a unique symmetric $\sigma$-derivation $\delta_{\sigma}$ such that $\delta_{\sigma}(x)=1$.
With the notations of example \ref{xp4} above, it is given by $\delta_{\sigma}(x^n) = [n]_{v} x^{n-1}$ with $[n]_{v} := (n)_{v^2}/v^{n-1}$.
If we define $\tau$ by the property $\tau(x) = vx$, one easily sees that $\partial_{\sigma} = \tau \delta_{\sigma}$ and $\delta_{\sigma} = \tau^{-1} \partial_{\sigma}$.
\end{enumerate}
\end{rmks}

The asymmetry in the $\sigma$-derivation formula has strong consequences:

%%%%%%%%%%
\begin{lem} \label{disym}
Let $\sigma$ be an endomorphism of a commutative $R$-algebra $A$ and $D : A \to M$ a $\sigma$-derivation.
\begin{enumerate}
\item Assume that there exists $x \in A$ such that $D(x)$ is regular in $M$ (not a torsion element), then
\begin{displaymath}
\forall z \in A, \quad D(z) = 0 \Rightarrow \sigma(z) = z.
\end{displaymath}
\item Assume that there exists $x \in A$ such that $y := x -\sigma(x)$ is regular on $A$ (not a zero divisor and not zero either), then
\begin{displaymath}
\forall z \in A, \quad \sigma(z) = z \Rightarrow D(z) = 0.
\end{displaymath}
\end{enumerate}
\end{lem}

\begin{pf}
We have for all $x, z \in A$,
\begin{displaymath}
xD(z) + \sigma(z)D(x) = D(xz) = D(zx) = zD(x) + \sigma(x)D(z).
\end{displaymath}
Thus, if we set $y := x - \sigma(x)$, we see that $yD(z) = (z - \sigma(z))D(x)$ and both assertions follows.
$\quad \Box$
\end{pf}

If $R \to R'$ is any homomorphism, then a $\sigma$-derivation $A \to M$ over $R$ extends by linearity to a $1 \otimes \sigma$-derivation $R' \otimes_{R} A \to R' \otimes_{R} M$ over $R'$.
Conversely if we are given some $R_{0} \to R$, then any $\sigma$-derivation of $A$ over $R$ may be seen as a $\sigma$-derivation over $R_{0}$. Also, if $E' \to E$ is any injective  map, then any $\underline \sigma_{|E'}$-derivation gives rise to a $\underline \sigma$-derivation.

%%%%%%%%%%%%%%
\begin{prop} \label{locder}
If $A$ is an $E$-twisted $R$-algebra and $B$ is an $E$-twisted localization of $A$, then any $\sigma$-derivation $D : A \to M$ will extend uniquely to a $\sigma$-derivation (still denoted) $D : B \to B \otimes_{A} M$.
\end{prop}

\begin{pf}
We can write $B = S^{-1}A$ were $S$ is a submonoid of $A \setminus \{0\}$ such that $\sigma(S) \subset S$.
For $x \in S$, we must have
\begin{equation}
0 = D(1) = D\left(x  \frac 1x\right) = x D\left(\frac 1x\right) + \frac{1}{\sigma(x)} D(x)
\end{equation}
and consequently $\displaystyle D\left(\frac 1x\right) = - \frac {D(x)}{x\sigma(x)}$.
It follows that if $y \in A$, we will have
\begin{equation} \label{locfrm}
D\left(\frac yx\right) = y D\left(\frac 1x\right) + \frac{1}{\sigma(x)} D(y) = - y\frac {D(x)}{x\sigma(x)} + \frac{1}{\sigma(x)} D(y) = \frac {xD(y)-yD(x)}{x\sigma(x)}.
\end{equation}
This proves uniqueness, and existence can then be checked directly.
$\quad \Box$
\end{pf}

If $A$ is an $E$-twisted commutative $R$-algebra, we will denote by $\mathrm{Der}_{R,\underline \sigma}(A,M)$ the set of all $\underline \sigma$-derivations of $A$ into $M$ and let
\begin{displaymath}
\mathrm T_{A/R,\underline \sigma} := \mathrm{Der}_{R, \underline \sigma}(A,A).
\end{displaymath}
Note that if $R \to R'$ is any map and $A' := R' \otimes_{R} A$, there exists a canonical map
$$
\mathrm{Der}_{R,\underline \sigma}(A,M) \to \mathrm{Der}_{R',\underline \sigma}(A',A' \otimes_{A} M).
$$

Also, as usual, when $A$ is trivially twisted, we will drop the prefix $\underline \sigma$ everywhere.

%%%%%%%%%%
\begin{prop}
If $A$ is an $E$-twisted $R$-algebra and $M$ is an $A$-module, then $\mathrm{Der}_{R,\underline \sigma}(A,M)$ is an $A$-submodule of $\mathrm{Hom}_{R}(A,M)$ for the action on the left.
\end{prop}

\begin{pf}
We have by definition
\begin{displaymath}
\mathrm{Der}_{R,\underline \sigma}(A,M) = \sum_{i \in E} \mathrm{Der}_{R,\sigma_{i}}(A,M),
\end{displaymath}
and we may therefore assume that we are in the $1$-twisted case.
Then our assertion follows immediately from the left $A$-linearity of the maps $D \mapsto [D, x]_{\sigma}$.
If $D, D' \in \mathrm{Der}_{R,\sigma}(A,M)$ and $x \in A$, we have
\begin{displaymath}
[D + D', x]_{\sigma} = [D, x]_{\sigma} + [D', x]_{\sigma} = D(x) + D'(x) = (D + D')(x).
\end{displaymath}
Also, if $D \in \mathrm{Der}_{R,\sigma}(A,M)$ and $x, y \in A$, we have
\begin{displaymath}
[y \circ D, x]_{\sigma} = y \circ [D, x]_{\sigma} = y \circ D(x) = (y \circ D)(x). \quad \Box
\end{displaymath}
\end{pf}

%%%%%%%%%%%%%%%
\begin{rmks}
\begin{enumerate}
\item In general, $\mathrm{Der}_{R,\underline \sigma}(A,M)$ is \emph{not} an $A$-submodule for the action on the right.
\item
The construction of $\mathrm{Der}_{R, \underline \sigma}(A,M)$ is not functorial in $A$ and does not even commute with extensions of $R$ in general.
\end{enumerate}
\end{rmks}

%%%%%%%%%%%%%%%%%%%%%%
\begin{xmps}
\begin{enumerate}
\item If $\sigma$ is an endomorphism of $A := R[x]$ (or $A := R[x,x^{-1}]$, or $A := K(x)$ if $R = K$ is a field), then $\mathrm T_{A/R, \sigma}$ is a free module of rank one.
\item Assume that $R$ is an integral domain and that $A := R[x_{1}, \dots, x_{n}]$ is endowed with $n$ endomorphism $\sigma_{i}$ such that $\sigma_{i}(x_{j}) = x_{j}$ whenever $i \neq j$ and $\sigma(x_{i}) \neq x_{i}$ for all $i = 1 \ldots, n$.
Then, $\mathrm T_{A/R, \underline \sigma}$ is a free module of rank $n$.
More precisely, a basis is given by by the $\partial_{i}$'s of \eqref{sigformn}.
\end{enumerate}
\end{xmps}

In order to check the last example, one may use the following result:

%%%%%%%%%%%%%%%%
\begin{prop}
Let $\sigma$ be an endomorphism of an $R$-algebra $A$ and $R' := A^{\sigma=1}$.
Let $M$ be an $A$-module.
Assume that there exists $x \in A$ such that $y := x -\sigma(x)$ is regular on $M$.
Then, we have
\begin{displaymath}
\mathrm{Der}_{R,\sigma}(A,M) = \mathrm{Der}_{R',\sigma}(A,M)
\end{displaymath}
\end{prop}

\begin{pf}
Follows from lemma \ref{disym}.
$\quad \Box$
\end{pf}

Although the module of $\underline \sigma$-derivations does not have a very nice functorial behavior, we have the following:

%%%%%%%%%%%%%%
\begin{prop} \label{loctf}
If $A$ is an $E$-twisted $R$-algebra of finite type, $B$ is an $E$-twisted localization of $A$ and $M$ is an $A$-module, then there exists an isomorphism
\begin{displaymath}
B \otimes_{A} \mathrm{Der}_{R, \underline \sigma}(A,M) \simeq \mathrm{Der}_{R, \underline \sigma}(B, B \otimes_{A} M).
\end{displaymath}
\end{prop}

\begin{pf}
It follows from proposition \ref{locder} that there exists a map
\begin{displaymath}
B \otimes_{A} \mathrm{Der}_{R, \underline \sigma}(A,M) \to \mathrm{Der}_{R, \underline \sigma}(B, B \otimes_{A} M).
\end{displaymath}
which is easily seen to be $A$-linear and injective (since $B$ is flat over $A$).
Surjectivity will follow if we show that any $\underline \sigma$-derivation $D'$ into $B \otimes_{A} M$ has the form $\frac 1y D$ where $D$ is a $\underline \sigma$-derivation into $M$ and $y \in S$.
But if we choose a set of generators $\{x_{1}, \ldots, x_{n}\}$ of $A$, we will have $D'(x_{i}) \in \frac 1{y_{i}} M$ for some $y_{i} \in S$ and we can take $y := \prod_{i = 1}^n y_{i}$.
$\quad \Box$
\end{pf}

Recall that, if $A$ is a commutative $R$-algebra, then an \emph{$R$-linear action} of an $A$-module $T$ on an $A$-module $M$ is a biadditive map
\begin{displaymath}
T \times M \to M, \quad (D, s) \mapsto Ds.
\end{displaymath}
such that
\begin{displaymath}
\forall D \in T, \forall s \in M, \quad \forall x \in A, \quad (xD)s = x(Ds) \quad \mathrm{and} \quad \forall a \in R, \quad D(as) = a(Ds).
\end{displaymath}
We will usually denote by $D_{M} : M \to M$, the map defined by $D_{M}(s) = Ds$.

If both $M$ and $N$ are endowed with an $R$-linear action of $T$, then a \emph{$T$-horizontal map} is an $A$-linear map $u : M \to N$ such that
\begin{displaymath}
\forall D \in T, \forall m \in M, \quad Du(m) = u(Dm).
\end{displaymath}
The category obtained this way is equivalent to the category of (left) $U$-modules where $U$ denotes the tensor algebra of $T$ over $A$.
In particular, this is an abelian category with sufficiently many injective and projective objects.

For example, an $E$-twisted commutative $R$-algebra $A$ is naturally endowed with an action of $\mathrm T_{A/R, \underline \sigma}$.

%%%%%%%%%%%%%%
\begin{dfn}
Let $A$ be an $E$-twisted commutative $R$-algebra and $M$ an $A$-module endowed with an \emph{$R$-linear action} of $\mathrm T_{A/R, \underline \sigma}$.
Then, a \emph{horizontal section} of $M$ is an element $s$ such that
\begin{displaymath}
\forall D \in \mathrm T_{A/R,\underline \sigma}, \quad Ds = 0.
\end{displaymath}
\end{dfn}

Clearly, the set $M'$ of all horizontal sections of $M$ may be identified with the set of all horizontal maps $A \to M$.
Also, we have
\begin{displaymath}
M' = \bigcap_{D \in \mathrm T_{A/R, \underline \sigma}} \ker D_{M}.
\end{displaymath}
In particular, we see that $M'$ is a $R$-submodule of $M$.

%%%%%%%%%
\begin{prop}
If $A$ is an $E$-twisted commutative $R$-algebra, then the horizontal sections of $A$ form an $R$-subalgebra $A'$ of $A$ which is contained in the ring $A^{\underline \sigma = \underline 1}$ of invariants of $A$.
\end{prop}

\begin{pf}
The first assertion follows from the $\sigma$-Leibnitz rule \eqref{sigleib} and the second one from lemma \ref{easy} which tell us that for all $i \in E$, the map $1 -\sigma_{i}$ is a $\underline \sigma$-derivation of $A$.
$\quad \Box$
\end{pf}

%%%%%%%%%%%%%%%%%%%%%
\begin{xmps}
\begin{enumerate}
\item
Recall from \cite{LeStumQuiros15} that if $q \in R$, then the $q$-\emph{characteristic} of $R$ is the smallest positive integer $p =: q\mathrm{-char} (R)$ such that $1 + q + \cdots q^{p-1}= 0$ if it exists, and zero otherwise.
Also, $R$ is said to be \emph{$q$-flat} (resp. \emph{$q$-divisible}) if $1 + q + \cdots + q^{m-1}$ is regular (resp. invertible) in $R$ or equal to $0$.

Now, we let $\sigma$ be the endomorphism of $R[x]$ defined by $\sigma(x) = qx$ and we assume that $R$ is $q$-flat.
Then,
\begin{enumerate}
\item If $q\mathrm{-char}(R) = 0$, then the set $A'$ of horizontal sections of $A$ is exactly $R$.
\item If $q\mathrm{-char}(R) = p > 0$, we have $A' = R[x^p]$.
\end{enumerate}
Note that, in both cases, we will actually have $A' = A^{\sigma = 1}$ if we also assume that $1 -q$ is not a zero-divisor.
\item We let $\sigma$ be the endomorphism of $R[x]$ defined by $\sigma(x) = x + 1$.
\begin{enumerate}
\item If $\mathrm{Char}(R) = 0$, then the set $A'$ of horizontal sections of $A$ is exactly $R$.
\item If $\mathrm{Char}(R) = p > 0$, we have $A' = R[x^p -x]$.
\end{enumerate}
Again, in both cases, we have $A' = A^{\sigma = 1}$.
\end{enumerate}
\end{xmps}

%%%%%%%%%%%%%%
\begin{dfn} \label{Liebnitz}
Let $A$ be a commutative $R$-algebra, $\sigma$ a ring endomorphism of $A$ and $M$ an $A$-module.
An $R$-linear endomorphism $D_{M}$ of $M$ is called a \emph{$\sigma$-derivation of $M$} if there exists $D_{A} \in \mathrm T_{A/R,\sigma}$ such that
\begin{displaymath}
\forall x \in A, \quad [D_{M}, x]_{\sigma} = D_{A}(x).
\end{displaymath}
More generally, if $A$ is an $E$-twisted commutative $R$-algebra, then a \emph{$\underline \sigma$-derivation of an $A$-module $M$} is a finite sum $D = \sum_{i \in E} D_{i}$ of $\sigma_{i}$-derivations of $M$.
Finally, a \emph{twisted differential $A$-module} is an an $A$-module $M$ endowed with an \emph{$R$-linear action} of $\mathrm T_{A/R, \underline \sigma}$ by $\underline \sigma$-derivations of $M$.
\end{dfn}

Note that if $D_{M}$ is a \emph{$\sigma$-derivation of $M$}, then $D_{A}$ is uniquely determined by $D_{M}$.
When we want to emphasize the role of $D_{A}$, we will call $D_{M}$ a \emph{$D_{A}$-derivation}.
Note that, by definition, $D_{M}$ is a $D_{A}$-derivation if and only if it satisfies the \emph{twisted Leibnitz rule}:
\begin{equation} \label{leibmod}
\forall x \in A, \forall s \in M, \quad D_{M}(xs) = D_{A}(x)s + \sigma(x)D_{M}(s).
\end{equation}

It should also be remarked that in the case $M = A$, this is compatible with definition \ref{sigmadev}.

It sounds natural to define a twisted differential operator as an operator that can be built by standard operations from functions and twisted derivations.
However, as it is already the case in the untwisted case, we obtain a ring which is smaller than expected in general.

%%%%%%%%%%
\begin{dfn} \label{smdif}
Let $A$ be an $E$-twisted commutative $R$-algebra.
Then, the ring of \emph{small twisted differential operators} $\overline{\mathrm D}_{A/R,\underline \sigma}$ is the subring of $\mathrm{End}_{R}(A)$ generated by $A$ and $\mathrm T_{A/R, \underline \sigma}$.
\end{dfn}

%%%%%%%%%%%%%%%%
\begin{xmps}
\begin{enumerate} 
\item Assume that $R = \mathbb C$ and $A = \mathbb C[x]$.
\begin{enumerate}
\item If $\sigma = \mathrm{Id}_{A}$, then $\overline{\mathrm D}_{A/R,\sigma}$ is nothing but the usual ring of differential operators on the affine line (the Weyl algebra in one variable).
\item If $\sigma(x) = x+h$ where $0 \neq h \in \mathbb C$, then $\overline{\mathrm D}_{A/R,\sigma}$ is the usual ring of finite difference operators on $\mathbb C[x]$.
\item If $\sigma(x) = qx$ where $q \in \mathbb C$ satisfies $|q| \neq 0,1$, then $\overline{\mathrm D}_{A/R,\sigma}$ is the usual ring of $q$-difference operators on $\mathbb C[x]$ (the quantum Weyl algebra).
\end{enumerate}
\item There are however some situations where the ring of \emph{small} twisted differential operators deserves its name.
\begin{enumerate}
\item As it is well known, if $R = \mathbb F_{2}$, $A = \mathbb F_{2}[x]$ and $\sigma = \mathrm{Id}_{A}$, then $\overline{\mathrm D}_{A/R,\sigma} = A \oplus \mathrm T_{A/R,\sigma}$.
In particular, this is \emph{not} the Weyl algebra.
\item In the same situation $R = \mathbb F_{2}$, $A = \mathbb F_{2}[x]$ but $\sigma (x) = x+1$, we also get $\overline{\mathrm D}_{A/R,\sigma} = A \oplus \mathrm T_{A/R,\sigma}$.
\item We want to emphasize the fact that his phenomenon also occurs in characteristic zero: if $R = \mathbb C$, $A = \mathbb C[x]$ and $\sigma(x) = -x$, then we have $\overline{\mathrm D}_{A/R,\sigma} = A \oplus \mathrm T_{A/R,\sigma}$ again.
We understand that this is a consequence of the $q$-characteristic being positive.
\end{enumerate}
\item Assume that $R = \mathbb C[t,t^{-1}]$ and $A = R[x_{1}, \cdots, x_{n}]$ is endowed with the endomorphism given by $\sigma_{i}(x_{j}) = \delta_{ij}tx_{j}$ where $\delta$ denotes the Kronecker symbol.
Then $\overline{\mathrm D}_{A/R,\underline \sigma}$ is identical to the $n$-th quantized Weyl algebra denoted by $A^{(t)}_{n}$ defined by Backelin on page 319 of \cite{Backelin11}.
\item Assume that $R=K$ is a field containing $\mathbb Q(q)$ with $q \not \in \mathbb Q^{\mathrm{alg}}$ and let $A = K[x]$.
We endow the set $E = \{0, 1,-1\}$ with the conditions $T_{0} = 1$ and $T_{1}T_{-1} = T_{-1}T_{1} = 1$,  and the ring $A$ with the endomorphisms given by $\sigma_{0}(x) = x$, $\sigma_{1}(x) = qx$ and $\sigma_{-1}(x) = \frac 1q x$.
Then $\overline{\mathrm D}_{A/R,\underline \sigma}$ is exactly the ring $D_{q}$ described in the article \cite{IyerMcCune02} (they use the letter $R$ for what we call $A$).
\item Let $\sigma$ be an endomorphism of an $R$-algebra $A$ and $\{\sigma_{n}\}_{n \in S}$ be a (compatible) system of roots of $\sigma$.
Then, we do \emph{not} have in general
\begin{displaymath}
\overline{\mathrm D}_{A/R,\underline \sigma} = \cup_{n \in S}\overline{\mathrm D}_{A/R,\sigma_{n}}.
\end{displaymath}
It works fine in the case $R = \mathbb Q (t)$, $A = R[x]$, $S = \{1,2\}$, $\sigma(x) = t^2x$ and $\sigma_{2}(x) = tx$ because we will have
\begin{displaymath}
\partial_{\sigma} = \partial_{\sigma_{2}} + \frac {t(t-1)}{t+1} x \partial^2_{\sigma_{2}}
\end{displaymath}
as one can easily check.
But the same formula shows that the result will not hold, for example, if we restrict to $R = \mathbb Q [t]$, because of the existence of a denominator.
\end{enumerate}
\end{xmps}

%%%%%%%%%%%%%%%%%%%
\begin{prop}
let $A$ be an $E$-twisted commutative $R$-algebra.
If $M$ is a $\overline{\mathrm D}_{A/R,\underline \sigma}$-module, then the induced action of $\mathrm T_{A/R,\underline \sigma}$ on $M$ is an action by $\underline \sigma$-derivations.
\end{prop}

\begin{pf}
This follows immediately from the definitions.
$\quad \Box$
\end{pf}

\begin{rmks}
\begin{enumerate}
\item
The converse of the assertion of the proposition is false in general since there already exists a problem in the trivially twisted situation in positive characteristic ($p$-curvature phenomenon).
\item It also is important to notice that $\overline{\mathrm D}_{A/R,\sigma}$ does not commute with extensions of $R$ in general, as it is already the case in the trivially twisted situation also.
\item It follows from proposition \ref{loctf} that if $A$ is an $E$-twisted $R$-algebra of finite type and $B$ is an $E$-twisted localization of $A$, then there exists an isomorphism
\begin{displaymath}
B \otimes_{A} \overline{\mathrm D}_{A/R,\underline \sigma} \simeq \overline{\mathrm D}_{B/R,\underline \sigma}.
\end{displaymath}
\end{enumerate}
\end{rmks}

%%%%%%%%%%%%%%%%%%%%%%%%%%%%%
%%%%%%%%%%%%%%%%%%%%%%%%%%%%%%
\section{Twisted differential algebras}

We introduce in this section the notion of twisted differential algebra.
This is a twisted algebra enhanced with a family of twisted derivations.

%%%%%%%%%%%%%%
\begin{dfn} \label{sigdiff}

An \emph{$E$-twisted differential $R$-algebra} is an $E$-twisted commutative $R$-algebra $(A, \underline \sigma_{A})$ endowed with a family of $\sigma_{A,i}$-derivations $D_{A,i}$ for all $i \in E$.
It is said of \emph{Schwarz type} if 
\begin{equation} \label{schcon}
\forall i, j \in E, \quad D_{A,i} \circ D_{A,j} = D_{A,j} \circ  D_{A,i} \quad  \mathrm{and} \quad   i \neq j \Rightarrow \sigma_{A,i} \circ D_{A,j} = D_{A,j} \circ \sigma_{A,i}.
\end{equation}
A \emph{morphism} between two $E$-twisted differential $R$-algebras $A$ and $B$ is a morphism $\varphi : A \to B$ of $E$-twisted $R$-algebras such that
\begin{displaymath}
\forall i \in E, \forall x \in A, \quad \varphi(D_{A,i}(x)) = D_{B,i}(\varphi(x)).
\end{displaymath}
\end{dfn}

Clearly, we obtain a category.
Although a twisted differential $R$-algebra is a triple $(A, \underline \sigma_{A}, \underline D_{A})$, we will normally use  only the first letter $A$ to denote it.
Also we will drop the index $A$ when there is no risk of confusion and write $D_{i}$ instead of $D_{A,i}$.
Note, and this is important, that this notion strongly depends on the choice of $E$ and \emph{not} only on $G(E)$.

%%%%%%%%%%%%%%%%%%%%%%
\begin{xmps}
\begin{enumerate}
\item If the $R$-algebra $A := R[x_{1}, \ldots, x_{n}]$ is endowed with $n$ endomorphisms $\sigma_{i}$ satisfying $\sigma_{i}(x_{j}) = x_{j}$ for $j \neq i$, we may use the $\sigma_{i}$-derivations $\partial_{i}$ of \eqref{sigformn}.
This defines an $n$-twisted differential algebra of Schwarz type.
\item If $R = K$ is a field and $A := K(x)$ is endowed with an endomorphism $\sigma$, then we may use the unique $\sigma$-derivation $\partial_{\sigma}$ such that $\partial_{\sigma}(x) = 1$ in order to turn $A$ into a $1$-twisted differential algebra.
This may be extended to several variables.
\item If $A := R[x,x^{-1}]$ and $\sigma(x) = qx$ where $q \in R^\times$ then we may use the $\sigma$-derivation $\partial_{\sigma}$ given by $\partial_{\sigma}(x^n) = (n)_{q}x^{n-1} $.
Again, it works as well with several variables.
\item Let $S$ be a set of positive integers (with $1 \in S$, as usual), $A= R[x]$, $\{\sigma_{n}\}_{n \in S}$ a compatible system of roots of some $R$-endomorphism $\sigma$ of an $R$-algebra $A$ and for each $n \in S$, $\partial_{n}$ the unique $\sigma_{n}$-derivation of $A$ such that $\partial_{n}(x) = 1$.
Then, the $S$-twisted $R$-algebra $(A, \underline \sigma, \underline D)$ is not of Schwarz type (at all) in general.
\end{enumerate}
\end{xmps}

%%%%%%%%%%%%%%%%%%%%%%
\begin{rmk}
The category of $E$-twisted differential $R$-algebras may be seen as a full subcategory of the category of (generalized) differential rings of Andr\'e (\cite{Andre01}, 2.1.2.1).
More precisely, if we denote by $A_{i}$ the $\sigma_{i}$-sesquimodule deduced from $A$, then
\begin{displaymath}
\underline D_{A} : A \to \prod_{i \in E} A_{i}
\end{displaymath}
is a generalized differential ring in Andr\'e's sense.
\end{rmk}

There exists a base change mechanism for twisted differential algebras: if $A$ is an $E$-twisted differential $R$-algebra (of Schwarz type) and $R \to R'$ is any homomorphism, then $A' := R' \otimes_{R} A$ has a natural structure of $E$-twisted differential $R'$-algebra (of Schwarz type).
And conversely, there exists an obvious restriction functor.
Also, if $E' \to E$ is any equivariant map, then there exists an obvious restriction functor from $E$-twisted differential $R$-algebras (of Schwarz type) to $E'$-twisted differential algebras (of Schwarz type).
Finally, any $E$-twisted localization of an $E$-twisted differential algebra is naturally an  $E$-twisted differential algebra.

%%%%%%
\begin{prop} Assume $E=E^{ab}$ (that is,  the conditions on $E$ include all the commutation conditions).
If $A$ is an $E$-twisted differential algebra of Schwarz type, then any $E$-twisted localization $B$ of $A$ is also of Schwarz type.
\end{prop}

\begin{pf}
In order to lighten the notations, we do no mention the indexes $A$ and $B$ in the $\sigma_{i}$'s.
We write $B = S^{-1}A$ where $S$ is a submonoid of $A \setminus \{0\}$ and we let $x \in S$ and $y \in A$.
Then, for $i \neq j$, we have, thanks to formula \eqref{locfrm}:
\begin{displaymath}
\sigma_{i} \left( D_{j}\left(\frac yx\right)\right) = \sigma_{i} \left( \frac {xD_{j}(y)-yD_{j}(x)}{x\sigma_{j}(x)} \right) = \frac { \sigma_{i}(x) \sigma_{i}(D_{j}(y)) -  \sigma_{i}(y) \sigma_{i}(D_{j}(x))}{ \sigma_{i}(x) \sigma_{i}(\sigma_{j}(x))}
\end{displaymath}
and
\begin{displaymath}
D_{j}\left(\sigma_{i}(\frac yx)\right) = D_{j}\left(\frac {\sigma_{i}(y)}{\sigma_{i}(x)}\right) = \frac {\sigma_{i}(x)D_{j}(\sigma_{i}(y))-\sigma_{i}(y)D_{j}(\sigma_{i}(x))}{\sigma_{i}(x)\sigma_{j}(\sigma_{i}(x))}.
\end{displaymath}
Thus, we see that the second condition of \eqref{schcon} is satisfied in $B$.
The first condition is shown to hold exactly in the same way.
$\quad \Box$
\end{pf}

%%%%%%%%%%
\begin{dfn} \label{smdif}
The module of \emph{derivations} of an $E$-twisted differential $R$-algebra $(A, \underline \sigma, \underline D)$ is the $A$-submodule $\mathrm T_{A/R, \underline \sigma, \underline D}$ of $\mathrm{End}_{R}(A)$ generated by the family $\underline D$. 
Given such an $(A, \underline \sigma, \underline D)$,  we will denote by $\overline{\mathrm D}_{A/R, \underline \sigma, \underline D}$  the subring  of $\mathrm{End}_{R}(A)$ generated by $A$ and $\mathrm T_{A/R, \underline \sigma, \underline D}$.
Finally, the \emph{module of differential forms} of $(A, \underline \sigma, \underline D)$ is the $A$-submodule $\Omega^1_{A/R, \underline \sigma, \underline D}$ of $A^E$ generated by the image of $\underline D : A \to A^E$.
\end{dfn}

%%%%%%%%%%%%%%
\begin{prop}
Let $A$ be an $E$-twisted differential $R$-algebra.
If $E$ is finite and $\Omega^1_{A/R, \underline \sigma, \underline D}$ is projective, then there exists an isomorphism
\begin{displaymath}
\mathrm{Hom}_{A}(\Omega^1_{A/R, \underline \sigma, \underline D}, A) \simeq \mathrm T_{A/R, \underline \sigma, \underline D}.
\end{displaymath}
\end{prop}

\begin{pf}
Since we assume that $E$ is finite, the image of the natural map
\begin{displaymath}
\mathrm{Hom}_{A}(A^E, A) \to \mathrm {End}_{R}(A), u \mapsto u \circ \underline D,
\end{displaymath}
is exactly $\mathrm T_{A/R, \underline \sigma, \underline D}$.
Now, by definition, $u$ is in the kernel of this map if and only if its restriction to $\Omega^1_{A/R, \underline \sigma, \underline D}$ is zero.
Finally, since we assume that $\Omega^1_{A/R, \underline \sigma, \underline D}$ is projective, the restriction map
\begin{displaymath}
\mathrm{Hom}_{A}(A^E, A) \to \mathrm{Hom}_{A}(\Omega^1_{A/R, \underline \sigma, \underline D}, A)
\end{displaymath}
is surjective and we obtain the expected identification.
$\quad \Box$
\end{pf}

%%%%%%%%%%%%%%
\begin{dfn} \label{twdim}
Let  $A$ be an $E$-twisted differential $R$-algebra.
Then, a \emph{$\underline D_{A}$-module} is an $A$-module $M$ endowed with a family of $D_{A,i}$-derivations $D_{M,i}$.
It is called \emph{integrable} if the $D_{M,i}$ commute with each other.
A \emph{morphism of $\underline D_{A}$-modules} is an $A$-linear map $u : M \to N$ such that
\begin{equation}
\forall i \in E, \forall s \in M, \quad u(D_{M,i}(s)) = D_{N,i}(u(s)). \label{comD}
\end{equation}
\end{dfn}

Clearly, with this definition, $\underline D_{A}$-modules form a category $\underline D_{A}\mathrm{-Mod}$.
We will denote the subcategory of integrable modules as $\underline D_{A}\mathrm{-Mod}^{\mathrm{int}}$.
Note that  $A$ itself is an $E$-twisted differential $A$-module and that it is integrable when $A$ is of Schwarz type.

%%%%%%%%%%%%%%
\begin{rmks}
\begin{enumerate}
\item The notion of (integrable) $E$-twisted differential $A$-module corresponds exactly to the notion of generalized (integrable) connection of Andr\'e (definition 2.2.1 of \cite{Andre01}).
\item If $A$ is an $E$-twisted differential $R$-algebra and $M$ and $N$ are two $E$-twisted differential $A$-modules, then there exists in general no natural structure of $E$-twisted differential $A$-module on $M \otimes_{A} N$ and $\mathrm{Hom}_{A}(M, N)$.
\end{enumerate}
\end{rmks}

It will follow from proposition \ref{abag2} below that if $A$ is an $n$-twisted differential $R$-algebra of Schwarz type, then $\underline D_{A}\mathrm{-Mod}$ is an abelian category with enough injectives.
Then, if $M$ is a $\underline D_{A}$-module, one sets
\begin{displaymath}
\mathrm R\Gamma_{\underline D}(M) := \mathrm{RHom}_{\underline D_{A}\mathrm{-Mod}}(A, M) \quad \mathrm{and} \quad \mathrm H^i_{\underline D}(M) := \mathrm{Ext}^i_{\underline D_{A}\mathrm{-Mod}}(A, M).
\end{displaymath}

%%%%%%%%%%%%%%%%%%%%
\begin{xmps}
\begin{enumerate}
\item We have
\begin{displaymath}
\mathrm H_{\underline D}^0(M) \simeq \cap_{i \in E} \ker D_{M,i}.
\end{displaymath}
\item If $(A, \sigma, D)$ is a $1$-twisted differential $R$-algebra and $(M, D)$ is a $1$-twisted differential $A$-module, one can show that
\begin{displaymath}
\mathrm R\Gamma_{D}(M) \simeq \left[
\xymatrix@R0cm{ M \ar[r]^{D} & M}\right].
\end{displaymath}
This is a complex concentrated in degree $0$ and $1$ whose cohomology is given by
\begin{displaymath}
\mathrm H_{D}^0(M) \simeq  \ker D \quad \mathrm{and} \quad \mathrm H_{\sigma}^1(M) \simeq M/\mathrm{Im}\, D.
\end{displaymath}
\end{enumerate}
\end{xmps}

Note that, if $A$ is an $E$-twisted differential $A$-algebra, then $A^{\underline D = \underline 0} := \mathrm H_{\underline D}^{0}(A)$ is a subring of $A$ called the ring of \emph{$\underline D$-constants} of $A$.

%%%%%%%%%%%%%%%%%%%%%%%%%%%%%
%%%%%%%%%%%%%%%%%%%%%%%%%%%%%%
\section{Twisted Weyl algebras}

In this section, we introduce the notion of twisted Weyl algebra which is attached to a twisted differential algebra in the same way as a a twisted polynomial ring is attached to a twisted ring.

We first recall some basic facts on filtrations and graduations.
We call \emph{graded $R$-algebra} (resp. \emph{filtered $R$-algebra}) an $R$-algebra $B$ endowed with a graduation (resp. an increasing filtration) by $R$-submodules which is indexed by a monoid $G$ endowed with a preorder compatible with multiplication on both sides (for example, the natural preorder defined by the conditions $g\leq gh$ and  $h \leq gh$ for all $g,h\in G$).

More precisely, a structure of graded $R$-algebra is given by a family of $R$-submodules $\mathrm{Gr}^g$ such that
\begin{displaymath}
B = \oplus_{g \in G} \mathrm{Gr}^{g}
\end{displaymath}
with the extra property that
\begin{displaymath}
\forall \varphi \in \mathrm{Gr}^g, \psi \in \mathrm{Gr}^h, \quad \varphi \psi \in \mathrm{Gr}^{gh} \quad \mathrm{and} \quad 1 \in \mathrm{Gr}^1.
\end{displaymath}
Note that $A := \mathrm{Gr}^1$ is an $R$-subalgebra of $B$ and that $B$ is a $G$-graded $A$-module.

A structure of filtered $R$-algebra is given by a family of $R$-submodules $\mathrm{Fil}^g$ such that
\begin{displaymath}
\forall g \leq h, \quad \mathrm{Fil}^{g} \subset \mathrm{Fil}^h \quad\mathrm{and} \quad B = \cup_{g \in G}\mathrm{Fil}^{g},
\end{displaymath}
with the extra properties that
\begin{displaymath}
\forall \varphi \in \mathrm{Fil}^g, \psi \in \mathrm{Fil}^h, \quad \varphi \psi \in \mathrm{Fil}^{gh} \quad \mathrm{and} \quad 1 \in \mathrm{Fil}^1.
\end{displaymath}
Note that $A := \mathrm{Fil}^1$ is an $R$-subalgebra of $B$ and that $B$ is a $G$-filtered $A$-module.

Any graded $R$-algebra may be seen as a filtered $R$-algebra by setting $\mathrm{Fil}^g := \oplus_{h\leq g} \mathrm{Gr}^h$.

%%%%%%%%%%%%%%%%%%%%
\begin{prop} Let $E$ a set with conditions and $G := G(E)$.
\begin{enumerate}
\item If $A$ is an $E$-twisted $R$-algebra, then $A[E]_{\underline \sigma}$ is a $G$-graded $R$-algebra for the grading $A[E]_{\underline \sigma} = \oplus_{g \in G} Ag$ and $G$ is a submonoid of (the multiplicative monoid of) $A[E]_{\underline \sigma}$.
\item  Conversely, let $B$ be a $G$-graded $R$-algebra with grading $B = \oplus_{g \in G} Ag$ where $G$ is a submonoid of (the multiplication monoid of) $B$.
Then, there exists a unique structure of $E$-twisted $R$-algebra on $A$ such that $B \simeq A[E]_{\underline \sigma}$ (as a graded $R$-algebras).
\end{enumerate}
\end{prop}

\begin{pf}
By definition, if $A$ is an $E$-twisted $R$-algebra, then
\begin{displaymath}
A[E]_{\underline \sigma} = \oplus_{g \in G} Ag
\end{displaymath}
as an $A$-module, $G$ is a submonoid for multiplication, and if we are given $g, h \in G$ and $x, y \in A$, we have $(xg)(yh) = x\sigma_{g}(y)gh$ which shows that multiplication is compatible with the grading.

Conversely, assume $B = \oplus_{g \in G} Ag$ is a graded $R$-algebra where $G$ is a submonoid of $B$.
If $x \in \mathrm{Gr^1} = A$ and $g \in G$, we must have $gx \in \mathrm{Gr}^g$ which means that $gx = \sigma_{g}(x)g$ for some $\sigma_{g}(x) \in A$.
Distributivity implies that $\sigma_{g}$ is an additive map, associativity that it is multiplicative, and the unit property shows that $\sigma_{g}$ is actually a ring homomorphism.
Associativity again shows that this defines an action of $G$ on $A$ (or, equivalently, an $E$-twisted structure).
The universal property of $A[E]_{\underline \sigma}$ provides a ring homomorphism $A[E]_{\underline \sigma} \to B$ which is clearly bijective.
$\quad \Box$
\end{pf}

Now, we want to relax the grading condition and replace it by a filtration condition:

%%%%%%%%%%%%%%
\begin{dfn} \label{weyl} Let $E$ be a set with conditions and $G := G(E)$.
An \emph{$E$-twisted Weyl algebra} is a filtered $R$-algebra $B$ such that there exists a commutative $R$-algebra $A$ and a morphism of monoids $G \to B, g \mapsto \partial_{g}$ inducing an isomorphism $\oplus_{g \in G} Ag \simeq B$ of filtered $A$-modules.
We call $\{\partial_{i}\}_{i \in E}$ a \emph{set of generators} of $B$ over $A$ and $\{\partial_{g}\}_{g \in G}$ a \emph{Weyl basis} of $B$ over $A$.
\end{dfn}

This definition only depends on $G$ and not on $E$.
 Note however that giving the map $G \to B$ is equivalent to giving a set of generators $\{\partial_{i}\}_{i \in E}$ satisfying the $E$-conditions in $B$.

The next result is in some sense a generalization of proposition 1.4 of \cite{Bourbaki12} (see also theorem 1.7.1 of \cite{Kassel95}) and may possibly be deduced from it but we'd rather give a direct proof. 

%%%%%%%%%%%%%%
\begin{thm} \label{oreweyl}
Let $B$ be an $n$-twisted Weyl algebra over $R$ with generators $\partial_{1}, \ldots, \partial_{n}$.
Then, there exists a unique $n$-twisted differential $R$-algebra $(A, \underline \sigma, \underline D)$ of Schwarz type such that for all $i = 1, \ldots, n$, we have
\begin{equation} \label{comsig2}
\forall x \in A, \quad \partial_{i} x = D_{i}(x) + \sigma_{i}(x) \partial_{i}.
\end{equation}
Conversely, given an $n$-twisted differential $R$-algebra of Schwarz type $(A, \underline \sigma, \underline D)$, there exists, up to isomorphism, a unique $n$-twisted Weyl algebra $B$ over $R$ with generator $\partial_{1}, \ldots, \partial_{n}$ over $A$ satisfying \eqref{comsig2}
\end{thm}

\begin{pf}
Let us begin with the first assertion.
Since multiplication respects the filtration, we can write for each $i = 1, \ldots, n$, a relation \eqref{comsig2} for some maps $\sigma_{i}, D_{i} : A \to A$.
Moreover, distributivity in $B$ implies that both maps are additive.
Now, we have on one hand for $x,y \in A$,
\begin{displaymath}
(\partial_{i} x) y = (D_{i}(x) + \sigma_{i}(x) \partial_{i}) y = D_{i}(x)y + \sigma_{i}(x) \partial_{i} y
\end{displaymath}
\begin{displaymath}
= D_{i}(x)y + \sigma_{i}(x) (\sigma_{i}(y)\partial_{i} + D_{i}(y)) = D_{i}(x)y + \sigma_{i}(x)D_{i}(y) +  \sigma_{i}(x)\sigma_{i}(y)\partial_{i}
\end{displaymath}
and on the other hand
\begin{displaymath}
\partial_{i} (xy) = D_{i}(xy) + \sigma_{i}(xy) \partial_{i}.
\end{displaymath}
The multiplication law in $B$ being associative, we must have
\begin{displaymath}
D_{i}(xy) = D_{i}(x)y + \sigma_{i}(x)D_{i}(y) \quad \mathrm{and} \quad  \sigma_{i}(xy) = \sigma_{i}(x)\sigma_{i}(y).
\end{displaymath}
Thus we see that $\sigma_{i}$ is multiplicative and $D_{i}$ satisfies the $\sigma_{i}$-Leibnitz rule.
Also, since $B$ is an $R$-algebra, we must have 
\begin{displaymath}
\forall a \in R, \quad a \partial_{i} = \partial_{i} a = D_{i}(a) + \sigma_{i}(a) \partial.
\end{displaymath}
It follows that $\sigma_{i}(a) = a$ and $D_{i}(a) = 0$ whenever $a \in R$.
Then $R$-linearity follows for both maps.
Now, for any $i, j \in \{1, \ldots, n\}$ and $x \in A$, we have
\begin{displaymath}
\partial_{i}\partial_{j} x = \partial_{i}D_{j}(x) + \partial_{i}\sigma_{j}(x) \partial_{j} = D_{i}(D_{j}(x)) + \sigma_{i}(D_{j}(x)) \partial_{i} + D_{i}(\sigma_{j}(x)) \partial_{j} + \sigma_{i}(\sigma_{j}(x)) \partial_{i} \partial_{j}
\end{displaymath}
Since $\partial_{i}$ and $\partial_{j}$ commute, we must have
\begin{displaymath}
D_{i} \circ D_{j} = D_{j} \circ D_{i} \quad \mathrm{and} \quad \sigma_{i} \circ \sigma_{j} = \sigma_{j} \circ \sigma_{i}
\end{displaymath}
but also when $i \neq j$,
\begin{displaymath}
\sigma_{i} \circ D_{j} = D_{j} \circ \sigma_{i}.
\end{displaymath}
Thus we see that all the expected commutation rules are satisfied.

The converse requires more work even if the uniqueness of $B$ is clearly automatic.
In order to show the existence of $B$, it is actually sufficient to build a \emph{big} $R$-algebra $C$ where all the expected relations hold and define $B$ as the subalgebra generated by $A$ and $\partial_{1}, \ldots, \partial_{n}$.
It will then automatically follow that any element of $B$ can be written as a finite sum $\sum x_{\underline k} \partial^{\underline k}$ with $x_{\underline k} \in A$ (we use the standard multiindex notation), and it will only remain to check that this expression is unique.
In other words, we will have to verify that $\sum x_{\underline k} \partial^{\underline k} \neq 0$ unless all $x_{\underline k} = 0$ (this is what \emph{big} means here).

Let $M$ be the free $A$-module on $\{e_{\underline k}\}_{\underline k \in \mathbb N^n}$ and $C = \mathrm{End}_{R}(M)$.
We may clearly consider $A$ as a subring of $C$ and define $\partial_{i}$ as the unique $R$-linear endomorphism of $M$ such that
\begin{displaymath}
\forall x \in A, \forall \underline k \in \mathbb N^n, \quad \partial_{i}(x e_{\underline k}) = D_{i}(x)e_{\underline k} + \sigma_{i}(x)e_{\underline k + \underline 1_{i}},
\end{displaymath}
where we denote by $\underline 1_{i}$ the element of $\mathbb N^n$ that has a $1$ in position $i$ and $0$ elsewhere. Let us check that formula  \eqref{comsig2} holds in $C$. If $x,y \in A$, we have 
\begin{displaymath}
(\partial_{i} \circ x)(y e_{\underline k}) = \partial_{i}(xy e_{\underline k}) = D_{i}(xy)e_{\underline k} + \sigma_{i}(xy)e_{\underline k+\underline 1_{i}}= (D_{i}(x)y + \sigma_{i}(x)D_{i}(y))e_{\underline k} + \sigma_{i}(x)\sigma_{i}(y)e_{\underline k+\underline 1_{i}}
\end{displaymath}
\begin{displaymath}
 = D_{i}(x)ye_{\underline k} + \sigma_{i}(x)(D_{i}(y)e_{\underline  k} + \sigma_{i}(y)e_{\underline k+ \underline 1_{i}}) = D_{i}(x)ye_{\underline k} + \sigma_{i}(x)\partial_{i} (ye_{\underline k})
\end{displaymath}
and it follows that
\begin{displaymath}
\partial_{i} \circ x = D_{i}(x) + \sigma_{i}(x) \circ \partial_{i}.
\end{displaymath}
Also, for all $i,j \in \{1, \ldots, n\}$ and $x \in A$, we have
\begin{displaymath}
\partial_{i}(\partial_{j}(x e_{\underline k})) = \partial_{i}(D_{j}(x)e_{\underline k} + \sigma_{j}(x)e_{\underline k + \underline 1_{j}})
\end{displaymath}
\begin{displaymath}
= D_{i}(D_{j}(x))e_{\underline k} + \sigma_{i}(D_{j}(x))e_{\underline k + \underline 1_{i}} + D_{i}(\sigma_{j}(x))e_{\underline k + \underline 1_{j}} + \sigma_{i}(\sigma_{j}(x))e_{\underline k + \underline 1_{i} + \underline 1_{j}}.
\end{displaymath}
With our commutation assumptions, this is symmetric in $i$ and $j$ and it follows that $\partial_{i}$ and $\partial_{j}$ commute with each other.
Finally, we have for all $i = 1, \ldots, n$,
\begin{displaymath}
\partial_{i}(e_{\underline k}) = D_{i}(1)e_{\underline k} + \sigma_{i}(1)e_{\underline k + \underline 1_{i}} = e_{\underline k + \underline 1_{i}}
\end{displaymath}
and it follows that for all for all $\underline l \in \mathbb N^n$, we have
\begin{displaymath}
(\sum x_{\underline k} \partial^{\underline k})(e_{\underline l}) = \sum x_{\underline k} e_{\underline k + \underline l}.
\end{displaymath}
Thus, we see that the condition $\sum x_{\underline k} \partial^{\underline k} = 0$ implies that all $x_{\underline k} = 0$ for all $\underline k \in \mathbb N^n$.
$\quad \Box$
\end{pf}

%%%%%%%%%%%%%%%%%%%%%%%%
\begin{dfn}
If $(A, \underline \sigma, \underline D)$ is an $n$-twisted differential $R$-algebra of Schwarz type, the unique $n$-twisted Weyl algebra $\mathrm D_{A/R,\underline \sigma, \underline D}$ over $R$ with generator $\partial_{1}, \ldots, \partial_{n}$ over $A$ such that formulas \eqref{comsig2} hold, is called the \emph{Ore extension} of $A$ by $\underline \sigma$ and $\underline D$.
\end{dfn}

%%%%%%%%%%%%%%%%%%%%%%%
\begin{xmps}
\begin{enumerate}
\item  When $A = R[x_{1}, \ldots, x_{n}]$, all $\sigma_{i} = \mathrm{id}_{A}$ and $D_{i} = \partial/\partial x_{i}$, then $\mathrm D_{A/R,\underline \sigma, \underline D}$ is the usual Weyl algebra over $R$.
\item More generally, when $A$ is a smooth $R$-algebra with an \'etale coordinates $x_{1}, \ldots, x_{n}$, all $\sigma_{i} = \mathrm{id}_{A}$ and $D_{i} = \partial/\partial x_{i}$,
then $\mathrm D_{A/R,\underline \sigma, \underline D}$ is the ring of differential operators \emph{of level $0$} of $A/R$ (not Grothendieck's ring of differential operators) that appears in crystalline cohomology.
Actually, when $R$ is a $\mathbb Q$-algebra, it is identical to the usual ring of differential operators.
 \end{enumerate}
\end{xmps} 

\begin{rmks}
\begin{enumerate}
\item Assume that $A = R[x]$ is endowed with an $R$-endomorphism $\sigma$ and that $D = \partial_{\sigma}$ is given by formula  \eqref{sigform} above.
Then, $\mathrm{D}_{A/R,\sigma, D}$ is the non commutative polynomial ring in two variables $R[x,\partial]$ with the commutation condition
\begin{displaymath}
\partial x = \sigma(x)\partial + 1.
\end{displaymath}
In the particular case $R = \mathbb Z[q]$ and $\sigma(x) = qx$, this is the quantum Weyl algebra used by M. Gros and the first author in \cite{GrosLeStum13}, definition 2.
\item Assume that $R = \mathbb C$ and $A = \mathbb C[x_{1}, \cdots, x_{n}]$ is endowed with
\begin{displaymath}
\sigma_{i}(x_{j}) = \delta_{ij}qx_{j} \quad \mathrm{and} \quad D_{i}(x^{\underline k}) = (k_{i})_{q}x^{\underline k -\underline 1_{i}}
\end{displaymath}
 for some $q \neq 0$.
Then $\mathrm D_{A/R,\underline \sigma, \underline D}$ is identical to the $n$-th quantized Weyl algebra denoted by $A^{(q)}_{n}$ on page 319 in \cite{Backelin11}.
\item Assume $R=K$ is a field containing $\mathbb Q(t)$, $A = K[x]$ and $\sigma(x) = tx$.
Then the ring of quantum differential operators over $K[x]$ introduced by Lunts and Rosenberg in \cite{LuntsRosenberg97} and studied by Iyer and McCune in \cite{IyerMcCune02} is \emph{not} the Weyl algebra in our sense.
It is bigger in general since it always contains usual differential operators.
Our ring is actually the $\sharp$ version introduced in section 1.2 of \cite{LuntsRosenberg97}.
\end{enumerate}
\end{rmks}

The $n$-twisted Weyl algebra satisfies the following universal property:

%%%%%%%%%%%%%%%%%%%%%%%
\begin{prop}
Let $A$ be an $n$-twisted differential $R$-algebra of Schwarz type.
Given a morphism of $R$-algebras $\varphi : A \to B$, and $n$ commuting elements $y_{1}, \ldots, y_{n} \in B$ such that
\begin{equation} \label{cuniv}
\forall i = 1, \ldots, n, \forall x \in A, \quad y_{i} \varphi(x) = \varphi(D_{i}(x)) + \varphi(\sigma_{i}(x)) y_{i},
\end{equation}
there exists a unique morphism of $R$-algebras $\Phi : \mathrm D_{A/R, \underline \sigma, \underline D} \to B$ that extends $\varphi$ and sends $\partial_{i}$ to $y_{i}$ for all $i = 1, \ldots, n$.
\end{prop}

\begin{pf}
By definition, we must have $\Phi(\sum x_{\underline k} \partial^{\underline k}) := \sum \varphi (x_{\underline k}) y^{\underline k}$ and uniqueness follows.
It only remains to show that the map is multiplicative.
This follows easily by induction from condition \eqref{cuniv} since the $y_{i}$'s commute with each other.
$\quad \Box$
\end{pf}

%%%%%%%%%%%%
\begin{prop} \label{abag2}
Let $A$ be an $n$-twisted differential $R$-algebra of Schwarz type.
If $M$ is a $\mathrm D_{A/R, \underline \sigma, \underline D}$-module, then the map 
\begin{displaymath}
\xymatrix@R0cm{ M  \ar[rr]^{D_{M,i}} && M 
\\ s \ar@{|->}[rr] && \partial _{i}s}
\end{displaymath}
turns $M$ into an integrable $n$-twisted differential $A$-module and we obtain an equivalence (an isomorphism) of categories
\begin{displaymath}
\mathrm D_{A/R, \underline \sigma, \underline D}\mathrm{-Mod} \simeq \underline D_{A}\mathrm{-Mod}^{\mathrm{int}}.
\end{displaymath}
\end{prop}

\begin{pf} This result follows immediately from the universal property of $\mathrm D_{A/R,\underline \sigma, \underline D}$ applied to the canonical morphism $A \to B := \mathrm{End}_{R}(M)$ when $M$ is an $A$-module.
More precisely, a structure of $\mathrm D_{A/R,\underline \sigma,\underline D}$-module on $M$ is given by $n$ commuting $R$-endomorphism $D_{M,i}$ of $M$ that satisfy
\begin{displaymath}
\forall i = 1, \ldots, n, \forall x \in A, \forall s \in M, \quad D_{M,i}(xs) = D_{i}(x)s + \sigma_{i}(x)D_{M,i}(s).
\end{displaymath}
This is exactly what we want.
$\quad \Box$
\end{pf}

As a consequence of the proposition, we also see that, if $M$ is any $n$-twisted integrable differential $A$-module, then the canonical morphism of $R$-algebras $\mathrm D_{A/R,\underline \sigma, \underline D} \to \mathrm{End}_{R}(M)$ takes values into the $R$-subalgebra generated by $A$ and the $\underline D$-derivations of $M$.
If we apply these considerations to the case $M = A$, we obtain:

%%%%%%%%%%%%%%%%%%
\begin{cor} \label{barsur}
There exists a morphism of $R$-algebras
\begin{displaymath}
\mathrm D_{A/R, \underline \sigma, \underline D} \to \overline{\mathrm D}_{A/R, \underline \sigma}
\end{displaymath}
from the Weyl algebra  to the ring of small  twisted differential operators
whose image is exactly the subring $\overline{\mathrm D}_{A/R, \underline \sigma, \underline D}$ generated by $A$ and the derivations (see definition \ref{smdif}). $\quad \Box$
\end{cor}

%%%%%%%%%%%%%%%%%%%%
\begin{cor}
If $A$ is an $n$-twisted differential $R$-algebra of Schwarz type, then the category of integrable $n$-twisted differential $A$-modules is an abelian category with sufficiently many projective and injective objects.$\quad \Box$
\end{cor}

Unlike the ring of small twisted differential operators, which does not even commute with extensions of $R$, the twisted Weyl algebra commutes with both extensions of $R$ and $A$:

%%%%%%%%%%%%%%%%%
\begin{prop}
\begin{enumerate}
\item
Let $A$ be an $n$-twisted differential $R$-algebra of Schwarz type, $R \to R'$ be any homomorphism and $A' := R' \otimes_{R} A$.
Then there exists a canonical isomorphism
\begin{displaymath}
R' \otimes_{R}\mathrm D_{A/R,\underline \sigma, \underline D} \simeq \mathrm D_{A'/R',\underline \sigma,\underline D}.
\end{displaymath}
\item
If $A \to B$ is a morphism of $n$-twisted differential $R$-algebras of Schwarz type, then there exists a canonical isomorphism
\begin{displaymath}
B \otimes_{A}\mathrm D_{A/R,\underline \sigma,\underline D} \simeq \mathrm D_{B/R, \underline \sigma,\underline D}.
\end{displaymath}

\end{enumerate}
\end{prop}

\begin{pf}
In both cases, we can use the universal property of $\mathrm D_{A/R,\underline  \sigma,\underline  D}$ in order to build the map and then compare basis on both sides.
$\quad \Box$
\end{pf}

%%%%%%%%%%%%%%%%%%%%%%%%%%%
%%%%%%%%%%%%%%%%%%%%%%%%%%
\section{Twisted coordinates} \label{twco}

In this last section, we are interested in twisted differential algebras where the twisted derivations may be seen as twisted partial derivatives.

We denote by $(A, \underline \sigma)$ an $n$-twisted commutative $R$-algebra: it means that we are given $n$ commuting endomorphisms $\sigma_{1}, \ldots, \sigma_{n}$ of an $R$-algebra $A$.

%%%%%%%%%%
\begin{dfn} \label{genD}
We say that $x_{1}, \ldots, x_{n} \in A$ are \emph{$\underline \sigma$-coordinates} for $A/R$ if
\begin{enumerate}
\item there exists for each $i = 1, \ldots, n$, a unique $\sigma_{i}$-derivation $\partial_{i}$ of $A/R$ such that 
\begin{equation} \label{unid}
\forall j=1, \dots, n, \quad \partial_{i} (x_{j}) = \left\{ \begin{array}{cl} 1 & \mathrm{if}\ j = i \\ 0 & \mathrm{otherwise} \end{array}\right.
\end{equation}
\item for any $\underline \sigma$-derivation $D$ of $A/R$, we have
\begin{equation}
D = \sum_{i =1} ^n D(x_{i})  \partial_{i}.
\end{equation}
\end{enumerate}

Then, $(A, \underline \sigma,  \underline \partial)$ is called a \emph{standard $n$-twisted differential $R$-algebra}.
The $\underline \sigma$-coordinates are said to be \emph{of Schwarz type} if the twisted differential algebra $(A, \underline \sigma, \underline \partial)$ is of Schwarz type.
\end{dfn}

Alternatively, in the definition of $\underline \sigma$-coordinates, one may require that the evaluation map
\begin{equation} \label{evalm}
\xymatrix@R0cm{T_{A/R, \underline \sigma} \ar[r] & A^n
\\ D \ar@{|->}[r] & (D(x_{1}), \ldots, D(x_{n}))}
\end{equation}
is bijective and that, for each $i = 1, \ldots, n$, the unique $\underline \sigma$-derivation satisfying \eqref{unid} is actually a $\sigma_{i}$-derivation.

%%%%%%%%%%%%%
\begin{xmps}
\begin{enumerate}
\item If $A$ is a smooth $R$-algebra and $\sigma_{i} = \mathrm{Id}_{A}$ for $i = 1, \ldots, n$, then $x_{1}, \ldots, x_{n}$ are $\underline \sigma$-coordinates for the $n$-twisted $R$-algebra $A$ if and only they are \'etale coordinates for $A$ over $R$ (they define an \'etale map between the polynomial algebra and $A$).
\'Etale coordinates are automatically of Schwarz type thanks to the usual Schwarz integrability conditions.
\item Assume that $R$ is an integral domain and that $A := R[x_{1}, \dots, x_{n}]$ is endowed with $n$ endomorphisms $\sigma_{i}$ such that $\sigma_{i}(x_{j}) = x_{j}$ whenever $i \neq j$ and $\sigma_{i}(x_{i}) \neq x_{i}$ for all $i = 1 \ldots, n$.
Then one can check that $x_{1}, \ldots, x_{n}$ are $\underline \sigma$-coordinates of Schwarz type for $A$.
\end{enumerate}
\end{xmps}

Actually, by specializing the endomorphisms, one may recover finite difference algebras, $q$-difference algebras as well as usual differential algebras:

%%%%%%%%%%%%%%%%%
\begin{dfn} Let $x_{1}, \ldots, x_{n}$ be $\underline \sigma$-coordinates for $A$.
\begin{enumerate}
\item If for all $i = 1, \ldots, n$, we have $\sigma_{i}(x_{i}) = x_{i} + h_{i}$ with $h_{i} \in R \setminus \{0\}$ and $\sigma_{i}(x_{j}) = x_{j}$ for $i \neq j$, we call $A$ a \emph{finite difference algebra} and say \emph{finite difference module} instead of $n$-twisted differential module on $A$.
\item If for all $i = 1, \ldots, n$, we have  $\sigma_{i}(x_{i}) = q_{i}x_{i}$ with $q_{i} \in R \setminus \{1\}$and $\sigma_{i}(x_{j}) = x_{j}$ for $i \neq j$, we call $A$ a \emph{$q$-difference algebra} and say \emph{$\underline q$-difference module} instead of a $n$-twisted differential module on $A$.
\item If for all $i, j \in \{ 1, \ldots, n\}$, we have  $\sigma_{i}(x_{j}) = x_{j}$, we call $A$ a \emph{(usual) differential algebra}, we say \emph{(usual) differential module} instead $n$-twisted differential module on $A$.
\end{enumerate}
\end{dfn}

%%%%%%%%%%
\begin{rmks}
\begin{enumerate}
\item When $n = 1$, an $x \in A$ is a $\sigma$-coordinate for $A$ if and only if the evaluation map $T_{A/R, \sigma} \to  A$ is bijective, and $\partial$ will then denote the unique $\sigma$-derivation such that $\partial(x) = 1$.
\item In the finite difference situation, it is common to write $\Delta_{h}$ instead of $\partial$.
Formula \eqref{class} below gives back the classical formula for the finite difference operator (when $h \in R^\times$):
\begin{displaymath}
\Delta_h(f)(x) = \frac {f(x+h) - f(x)}{h}.
\end{displaymath}
\item In the $q$-difference situation, one usually writes $\delta_{q}$ instead of $\partial$, and formula \eqref{class} below gives back the classical formula for the $q$-difference operator (when $1 -q \in R^\times$ and $x \in A^\times$):
\begin{displaymath}
\delta_q(f)(x) = \frac {f(qx) - f(x)}{qx - x}.
\end{displaymath}
\item In the usual situation, we do write $\partial$ and we have the classical formulas (whenever it means something):
\begin{displaymath}
\partial(f)(x) = \lim_{h\to 0}\Delta_h(f)(x) =  \lim_{q\to 1}\delta_q(f)(x).
\end{displaymath}
\end{enumerate}
\end{rmks}

Many more examples come from the following result:

%%%%%%%%%%%%%%
\begin{prop}
If $B$ is an $n$-twisted localization of $A$ and $x_{1}, \ldots, x_{n}$ are $\underline \sigma$-coordinates for $A/R$, then their images in $B$ are also $\underline \sigma$-coordinates for $B/R$.
\end{prop}

\begin{pf}
It follows from proposition \ref{loctf} that $T_{B/R, \underline \sigma} \simeq B \otimes_{A} T_{A/R, \underline \sigma}$ and the evaluation map \eqref{evalm} therefore stays an isomorphism after extension to $B$ (and the $\partial_{i}$'s are $\sigma_{i}$-derivations of $B/R$).
$\quad \Box$
\end{pf}

The existence of the twisted partial derivative has some important consequences as we will see shortly.

%%%%%%%%%%
\begin{lem} \label{weak}
Let $x_{1}, \ldots, x_{n}$ be $\underline \sigma$-coordinates for $A$.
Then, if we set $y_{i} := x_{i} - \sigma_{i}(x_{i})$ for all $i = 1, \ldots, n$, we have
\begin{enumerate}
\item $\forall j = 1, \ldots, n, \quad  i \neq j \Rightarrow \left(\sigma_{i}(x_{j}) = x_{j} \ \mathrm{and} \ \sigma_{i}(y_{j}) = y_{j}\right)$.
\item if $\underline x$ is of Schwarz type, then $\forall j = 1, \ldots, n, \quad  i \neq j \Rightarrow \partial_{i}(y_{j}) = 0$.
\end{enumerate}
\end{lem}

\begin{pf}
It follows from lemma \ref{disym} that $\sigma_{i}(x_{j}) = x_{j}$ for $j \neq i$.
Moreover, since $\sigma_{i}$ and $\sigma_{j}$ commute with each other, we also have
\begin{displaymath}
\sigma_{i}(y_{j}) = \sigma_{i}(x_{j} - \sigma_{j}(x_{j})) = \sigma_{i}(x_{j}) - \sigma_{i}\sigma_{j}(x_{j}) = x_{j} - \sigma_{j}\sigma_{i}(x_{j}) = x_{j} - \sigma_{j}(x_{j}) = y_{j}.
\end{displaymath}
Finally, if $\underline x$ is of Schwarz type, then we have $\partial_{i}\sigma_{j} = \sigma_{j}\partial_{i}$ for $j \neq i$, and it follows that
\begin{displaymath}
\partial_{i}(y_{j}) = \partial_{i}(x_{j} - \sigma_{j}(x_{j})) = - \partial_{i}\sigma_{j}(x_{j}) = - \sigma_{j}\partial_{i}(x_{j}) = 0. \quad \Box
\end{displaymath}
\end{pf}

In a standard $n$-twisted differential algebra, we can recover the family of endomorphisms $\underline \sigma$ from the family of twisted derivations $\underline \partial$.
More precisely, we have:

%%%%%%%%%%
\begin{prop} \label{weak}
If $x_{1}, \ldots, x_{n}$ are $\underline \sigma$-coordinates for $A$ and we set $y_{i} := x_{i} - \sigma_{i}(x_{i})$ for all $i = 1, \ldots, n$, then we have
\begin{equation} \label{coooord}
\sigma_{i} = 1 - y_{i} \partial_{i}.
\end{equation}
\end{prop}

\begin{pf}
Since $1 - \sigma_{i}$ is a $\underline \sigma$-derivation thanks to lemma \ref{easy}, we have
\begin{equation} \label{sig}
1 - \sigma_{i} =\sum_{j=1}^n (x_{j} - \sigma_{i}(x_{j})) \partial_{j} = y_{i}\partial_{i}
\end{equation}
because $\sigma_{i}(x_{j}) = x_{j}$ for $j \neq i$.
$\quad \Box$
\end{pf}

%%%%%%%%%%%%%%
\begin{dfn} \label{stcoo}
The $\underline \sigma$-coordinates $x_{1}, \ldots, x_{n}$ are said to be \emph{strong} if for all $i = 1, \ldots, n$, we have $y_{i} := x_{i} - \sigma_{i}(x_{i}) \in A^\times$.
We will also say that the $n$-twisted $R$-algebra is \emph{strong}.
\end{dfn}

If this is the case, we can recover the twisted partial derivative $\partial_{i}$ from the endomorphism $\sigma_{i}$ by the formula
\begin{equation} \label{class}
\partial_{i} = \frac 1{y_{i}} (1 -\sigma_{i}).
\end{equation}

%%%%%%%%%%%%%%%%%%%%%%
\begin{xmps}
\begin{enumerate}
\item Assume that $n = 1$ and  $A = R[x]$.
Then $x$ is a strong $\sigma$-coordinate for $A$ if and only if $\sigma(x) = x + h$ with $h \in R^\times$.
\item Assume that $n = 1$ and $A = K(x)$, where $R = K$ is a field.
Then $x$ is a strong $\sigma$-coordinate for $A$ if and only if $\sigma \neq \mathrm{Id}_{A}$.
\item Assume that $n = 1$, $A = R[x,x^{-1}]$ and $\sigma(x) = qx$ with $q \in R^\times$.
Then $x$ is a strong $\sigma$-coordinate for $A$ if and only if $1 - q \in R^\times$.
This is the case for example if $R = K$ is a field and $q \neq 0, 1$ or more generally, if $R$ contains a field $K$ with $q \in K \setminus \{0, 1\}$.
\end{enumerate}
\end{xmps}

Recall that we defined for the $n$-twisted algebra $(A, \underline \sigma)$,
\begin{displaymath}
A^{\underline \sigma = \underline 1} := \{z \in A, \forall i = 1, \ldots, n, \sigma_{i}(z) = z\} = \cap_{i=1}^n \ker (1 - \sigma_{i}),
\end{displaymath}
and for the $n$-twisted differential algebra $(A, \underline \sigma,  \underline \partial)$,
\begin{displaymath}
A^{\underline \partial = \underline 0} := \{z \in A, \forall i = 1, \ldots, n, D_{i}(z) = 0\} = \cap_{i=1}^n\ker \partial_{i}.
\end{displaymath}

%%%%%%%%%%%%%%%%%%%%%
\begin{prop} If $A$ is a standard $n$-twisted differential algebra, we have $A^{\underline \partial =\underline 0} \subset A^{\underline \sigma= \underline 1}$ with equality when $\underline x$ is strong.
\end{prop}

\begin{pf}
Follows form formulas \eqref{coooord} and \eqref{class}.
$\quad \Box$
\end{pf}

Recall from Section  \ref{twder} that we can associate to the $n$-twisted algebra $(A, \underline \sigma)$ its module of $\underline \sigma$-derivations $\mathrm T_{A/R,\underline \sigma}$ and its ring of small  twisted differential operators $\overline{\mathrm D}_{A/R,\underline \sigma}$.
On the other hand, if we have $\underline \sigma$-coordinates, we may also consider the module of derivations $\mathrm T_{A/R, \underline \sigma, \underline \partial}$ and the ring of small differential operators $\overline{\mathrm D}_{A/R, \underline \sigma, \underline \partial}$ associated to the $n$-twisted differential algebra $(A, \underline \sigma,  \underline \partial)$  (see definition \ref{smdif}).
It should be clear that, in both cases, the latter, which is always contained in the former, is actually equal to it.
Note also that we have here $\Omega^1_{A, \underline \sigma, \underline \partial} = A^n$.
More important for us, we have the following:

%%%%%%%%%%%%%%%%
\begin{prop} \label{sigmad} If $A$ is a standard $n$-twisted differential $R$-algebra, then the category of $A$-modules endowed with an $R$-linear action by derivation of $\mathrm T_{A/R, \underline \sigma}$ is equivalent (and even isomorphic) to the category $\underline \partial_{A}\mathrm{-Mod}$.
\end{prop}

\begin{pf}
This follows from the fact that $\{\partial_{1}, \ldots, \partial_{n}\}$ is a basis of $\mathrm T_{A/R, \underline \sigma}$.
$\quad \Box$
\end{pf}

%%%%%%%%%%%%%%%
\begin{dfn}
If $A$ is a standard $n$-twisted differential algebra, we call an $A$-module with a linear action of $\mathrm T_{A/R, \underline \sigma}$ by derivation \emph{integrable} if the corresponding $\underline \partial_{A}$-module is integrable.
\end{dfn}

They form a category that we will denote by $\mathrm{MIC}(A, \underline \sigma)$.
Note that this is not an intrinsic definition since it seems to rely on the choice of the $\underline \sigma$-coordinates.

%%%%%%%%%%%%%%%%%%%%%%%%%%%%%
\begin{dfn} If $x_{1}, \ldots, x_{n}$ are $\underline \sigma$-coordinates of Schwarz type for $A/R$, then the \emph{standard twisted Weyl algebra} associated to $\underline x$ is the Ore extension $\mathrm D_{A/R, \underline \sigma, \underline \partial} $ of $A$ by $\underline \sigma$ and $\underline \partial$.
\end{dfn}

We will also denote by $\partial_{1}, \ldots, \partial_{n}$ the corresponding generators of $\mathrm D_{A/R, \underline \sigma, \underline \partial}$.
They should not be confused with the $\underline \sigma$-derivations that we also call $\partial_{1}, \ldots, \partial_{n}$.
If necessary, we will write $\partial_{A,i}$ when we need to emphasize that we are considering the $\sigma_{i}$-derivation of $A$ (and not one of the generators of the $n$-twisted Weyl algebra).
Also, if we write as usual $y_{i} := x_{i} - \sigma(x_{i})$, we will set below $\sigma_{i} := 1 -y_{i}\partial_{i} \in \mathrm D_{A/R, \underline \sigma, \underline \partial}$.
Then, we will insist on the notation $\sigma_{A,i}$ when we want to emphasize the fact that we are considering the endomorphism of $A$ (and not the element of the twisted Weyl algebra).
This being said, we see that, by definition, an element of the standard $n$-twisted Weyl algebra $\mathrm D_{A/R, \underline \sigma, \underline \partial} $ may be uniquely written as a finite sum $\sum \underline z_{k} \underline  \partial^{\underline k}$, and multiplication is characterized by the properties
\begin{displaymath}
\partial_{i} z =  \partial_{A,i}(z) + \sigma_{A,i}(z) \partial_{i} \quad  \mathrm{and} \quad \partial_{i}\partial_{j} = \partial_{j} \partial_{i}.
\end{displaymath}

%%%%%%%%%%%%%%%%%%%
\begin{prop}\label{gengen}
If $A$ is a standard $n$-twisted differential $R$-algebra of Schwarz type, then there exists a canonical surjective map
\begin{equation}
\xymatrix@R0cm{
D_{A/R, \underline \sigma,  \underline \partial} \ar@{->>}[r]& \overline D_{A/R,\underline \sigma} \\ \partial_{i} \ar@{|->}[r] & \partial_{A, i}.} \label{tobar}
\end{equation}
Moreover, if we set $y_{i} = x_{i} - \sigma(x_{i})$ and $\sigma_{i} := 1 - y_{i} \partial_{i} \in \mathrm D_{A/R,\underline \sigma}$, then the image of $\sigma_{i}$ is $\sigma_{A,i}$.
\end{prop}

\begin{pf}
Follows from corollary \ref{barsur} and proposition \ref{weak} where the equality is meant to hold in $\mathrm T_{A/R, \underline \sigma} \subset \overline{\mathrm D}_{A/R, \underline \sigma}$ .
$\quad \Box$
\end{pf}

%%%%%%%%%%%%
\begin{lem} \label{comds}
Let $x_{1}, \ldots, x_{n}$ be $\underline \sigma$-coordinates of Schwarz type for $A$.
If for all $i = 1, \ldots, n$, we let $y_{i} = x_{i} - \sigma_{i}(x_{i}) \in A$ and  $\sigma_{i} := 1 - y_{i} \partial_{i} \in \mathrm D_{A/R, \underline \sigma, \underline \partial}$, then we have
\begin{enumerate}
\item $\forall z \in A, \quad \quad \sigma_{i} z = \sigma_{i}(z)\sigma_{i}$
\item $\forall j = 1, \ldots, n, \quad  i \neq j \Rightarrow \sigma_{i}\partial_{j} = \partial_{j}\sigma_{i}$
\item $\forall j = 1, \ldots, n, \quad  \sigma_{i}\sigma_{j} = \sigma_{j}\sigma_{i}$
\end{enumerate}
\end{lem}

\begin{pf}
The first equality is easily checked:
\begin{displaymath}
\sigma_{i}z = z - y_{i}\partial_{i}z = z - y_{i}\partial_{i}(z) + \sigma_{i}(z)y_{i}\partial_{i}
= \sigma_{i}(z) + \sigma_{i}(z)(1 -\sigma_{i}) = \sigma_{i}(z)\sigma_{i}.
\end{displaymath}
In order to prove the second identity, we compute
\begin{displaymath}
\sigma_{i}\partial_{j} = (1 - y_{i}\partial_{i})\partial_{j} = \partial_{j} - y_{i}\partial_{i}\partial_{j},
\end{displaymath}
and
\begin{displaymath}
\partial_{j} \sigma_{i} = \partial_{j}(1 - y_{i}\partial_{i}) = \partial_{j} - \partial_{j}y_{i}\partial_{i} = \partial_{j} - \partial_{j}(y_{i})\partial_{i} - y_{i}\partial_{j}\partial_{i} =  \partial_{j} - y_{i}\partial_{j}\partial_{i}.
\end{displaymath}
The last identity then easily follows:
the case $i = j$ is trivial and when $j \neq i$, we have
\begin{displaymath}
\sigma_{i}\sigma_{j} = \sigma_{i}(1 - y_{j}\partial_{j}) =  \sigma_{i} - \sigma_{i} y_{j}\partial_{j} = \sigma_{i} -  \sigma_{i} (y_{j})\sigma_{i}\partial_{j}
\end{displaymath}
\begin{displaymath} = \sigma_{i} -  y_{j}\sigma_{i}\partial_{j} = \sigma_{i} -  y_{j}\partial_{j}\sigma_{i} = (1 -  y_{j}\partial_{j})\sigma_{i} = \sigma_{j} \sigma_{i}. \quad \Box
\end{displaymath}
\end{pf}

%%%%%%%%%%%%
\begin{thm} \label{sigdif}
Let $(A, \underline \sigma)$ be an $n$-twisted differential $R$-algebra.
If $x_{1}, \ldots, x_{n}$ are $\underline \sigma$-coordinates of Schwarz type for $A$, then there exists a unique $A$-linear homomorphism of $R$-algebras
\begin{equation} \label{under1}
\xymatrix@R0cm{A[\underline T]_{\underline \sigma} \ar[r] & \mathrm D_{A/R, \underline \sigma, \underline \partial}
\\ T_{i} \ar@{|->}[r] & \sigma_{i} := 1 - (x_{i} - \sigma_{i}(x_{i})) \partial_{i}}
\end{equation}
If the $\sigma$-coordinates are strong, this map is an isomorphism.
\end{thm}

\begin{pf}
Both uniqueness and existence of the homomorphism will follow from the universal property of $A[\underline T]_{\underline \sigma}$ once we have checked that
\begin{displaymath}
\forall i, j \in \{1, \ldots, n\}, \quad \sigma_{i}\sigma_{j} = \sigma_{j}\sigma_{i}
\end{displaymath}
and 
\begin{displaymath}
\forall i = 1, \ldots, n, \forall z \in A, \quad \sigma_{i} z = \sigma_{i}(z)\sigma_{i}.
\end{displaymath}
This was shown in lemma \ref{comds}.
When $(A, \underline \sigma)$ is strong, one can show in the same way that there exists a unique $A$-linear homomorphism of $R$-algebras 
\begin{equation} \label{under2}
\xymatrix@R0cm{ \mathrm D_{A/R, \underline \sigma, \underline \partial} \ar[r] & A[\underline T]_{\underline \sigma}
\\ \partial_{i} \ar@{|->}[r] & \frac 1{y_{i}}(1 - T_{i})
}
\end{equation}
with $y_{i} := x_{i} - \sigma_{i}(x_{i})$.
This time, we use the universal property of $ \mathrm D_{A/R, \underline \sigma, \underline \partial}$.
We have to check that
\begin{equation} \label{comT}
\forall i, j \in \{1, \ldots, n\}, \quad \frac 1{y_{i}} (1 -T_{i})\frac 1{y_{j}} (1 -T_{j}) = \frac 1{y_{j}} (1 -T_{j}) \frac 1{y_{i}} (1 -T_{i})
\end{equation}
and
\begin{equation} \label{revm}
\forall i = 1, \ldots, n, \quad \forall z \in A, \quad \frac 1y_{i}(1 -T_{i}) z = \partial_{A,i}(z) + \sigma_{A,i}(z) \frac 1{y_{i}}(1-T_{i}).
\end{equation}
The identity \eqref{comT} is trivial when $j = i$.
When $j \neq i$, it follows from lemma \ref{weak} that $\sigma_{i}(\frac 1{y_{j}}) = \frac 1{y_{j}}$.
Therefore, $\frac 1{y_{j}}$ commutes with $T_{i}$ and the left hand side is symmetric on $i$ and $j$.
In order to show the identity \eqref{revm}, we first use formula \eqref{class} which tells us that $\partial_{A,i}(z) = \frac 1{y_{i}} (z - \sigma_{A,i}(z))$, and the condition can be rewritten as
\begin{displaymath}
(1 -T_{i})z = z - \sigma_{A,i}(z) + \sigma_{A,i}(z) (1 -T_{i}),
\end{displaymath}
which in turn reduces to $T_{i}z = \sigma_{i}(z)T_{i}$.
Finally, using the universal properties of both $R$-algebras, it is clear that the morphinisms of \eqref{under1} and \eqref{under2} are inverse to each other.
 $\quad \Box$
\end{pf}
 
%%%%%%%%%%%%%%%%%%%
 \begin{rmks}
 \begin{enumerate}
 \item  The homomorphism \eqref{under1} is actually a homomorphism of filtered $R$-algebras.
 \item If we compose the homomorphism \eqref{under1} with the surjection $\mathrm D_{A/R, \underline \sigma, \underline \partial} \to \overline{\mathrm D}_{A/R,\underline \sigma}$ and the inclusion $\overline{\mathrm D}_{A/R, \underline \sigma} \hookrightarrow \mathrm{End}_{R}(A)$, we obtain the canonical map
\begin{displaymath}
\xymatrix@R0cm{A[\underline T]_{\underline \sigma} \ar[r] & \mathrm{End}_{R}(A)
\\ T_{i} \ar@{|->}[r] & \sigma_{A,i}.}
\end{displaymath}
 \end{enumerate}
 \end{rmks}

%%%%%%%%%%%%%%%%%%%%%%
\begin{cor} \label{corr}
Let $(A, \underline \sigma)$ be an $n$-twisted differential $R$-algebra and $x_{1}, \ldots, x_{n}$ $\underline \sigma$-coordinates of Schwarz type for $A$.
If an $A$-module $M$ is endowed with an integrable action of $\mathrm T_{A, \underline\sigma}$, then for all $i = 1, \ldots, n$, the endomorphism defined by
\begin{displaymath}
\forall s \in M, \quad \sigma_{M,i}(s) =  s - (x_{i} -\sigma_{A,i}(x_{i})) \partial_{M,i}(s),
\end{displaymath}
is $\sigma_{A,i}$-linear and we obtain a functor
\begin{displaymath}
\mathrm{MIC}(A, \underline \sigma) \to \underline \sigma_{A}\mathrm{-Mod}.
\end{displaymath}
This is an equivalence of categories when $A$ is strong.
\end{cor}

Of course, when $A$ is strong, if we are given a $\underline \sigma_{A}$-module $M$, then the corresponding action will be given by
\begin{displaymath}
\forall s \in M, \quad \partial_{M,i}(s) = \frac {s - \sigma_{M,i}(s)}{x - \sigma_{A,i}(x)}.
\end{displaymath}

\begin{pf}
As we saw in proposition \ref{eqtwp}, there exists an equivalence of categories between $A[\underline T]_{\underline \sigma}$-modules and $n$-twisted $A$-modules, which in turn is equivalent to the category of $\underline \sigma_{A}\mathrm{-Mod}$.
We also showed in proposition \ref{abag2} that there exists an equivalence of categories between $\mathrm D_{A/R, \underline \sigma, \underline \partial}$-modules and $\underline \partial_{A}\mathrm{-Mod}^{\mathrm{int}}$, which by definition is equivalent to $\mathrm{MIC}(A, \underline \sigma)$.
The assertion therefore follows from theorem \ref{sigdif}. 
$\quad \Box$
\end{pf}
 
\begin{rmk}
Using this dictionary, we see that the modules considered by Claude Sabbah in \cite{Sabbah93} correspond to integrable twisted differential modules on $A := K[x_{1}, \ldots, x_{n}]$ with respect to $\partial_{A,i}(\underline x^{\underline k}) = q^{k_{i}} \underline x^{\underline k - 1_{i}}$.
\end{rmk}

Note also that we will always have a canonical map $\mathrm R\Gamma_{\underline D}M \to \mathrm R\Gamma_{\underline \sigma}M$ and that it is an isomorphism when $M$ is strong.

%%%%%%%%%%%%%%%%%%%%
\begin{xmp}
We can illustrate this equivalence of categories with the following classical case:
\begin{displaymath}
A = R[x,x^{-1}] \quad \mathrm{and}\quad  \sigma(x) = qx \ \mathrm{with}\ q, 1-q \in R^\times.
\end{displaymath}
Thus, we have
\begin{displaymath}
\forall n \in \mathbb Z, \quad \sigma(x^n) = q^nx^n \quad \mathrm{and}\quad  \partial(x^n) = (n)_{q} x^{n-1}.
\end{displaymath}
Then, there exists an isomorphism
\begin{displaymath}
\xymatrix@R0cm{ R[x,x^{-1}][T]_{\sigma} \ar[r]^-\simeq & \mathrm D_{R[x,x^{-1}]/R,\sigma, \partial}
\\ T \ar@{|->}[r] & 1 - (1 -q)x \partial \\ \frac {1 -T}{(1-q)x} & \partial \ar[l]}
\end{displaymath}
inducing an equivalence of categories
\begin{displaymath}
\left\{1\mathrm{-Twisted}\ \mathrm{differential}\  \mathrm{modules}\ \mathrm{on} \ R[x,x^{-1}] \right\} \simeq \left\{1\mathrm{-Twisted}\ \mathrm{modules}\ \mathrm{on} \ R[x,x^{-1}] \right\}.
\end{displaymath}
It is given by
\begin{displaymath}
\sigma_{M}(s) =  s  - (1-q)x \partial_{M}(s)
\quad
\mathrm{and}
\quad
\partial_{M,\sigma}(s) = \frac {s - \sigma_{M}(s)}{(1-q)x} .
\end{displaymath}
\end{xmp}
We want to stress out the fact that we do not make any assumption on $R$ which may as well be a ring of characteristic $p > 0$ nor on $q$ which may also be a $p$-th root of unity.
And of course, this example extends to higher dimension.
 
%\end{comment}

%%%%%%%%%%%%%%%%%%%%%%%%%%%
\bibliographystyle{plain}
\addcontentsline{toc}{section}{References}
\bibliography{BiblioBLS}

\end{document}